\documentclass[preprint,sort&compress,12pt]{elsarticle}
\usepackage{bbm}
\usepackage{amstext}
\usepackage{amsmath}
\usepackage{amssymb}
\usepackage{mathrsfs}
\usepackage{stmaryrd}
\usepackage{bm}
\usepackage{pstricks}
\usepackage{pst-coil}
\usepackage{pst-3d}
\usepackage{amsfonts}
\usepackage{amsthm} 
\newcommand{\mm}{\mathrm}
\newcommand{\ml}{\mathcal}

\newcommand{\be}{\begin{equation}}
\newcommand{\bea}{\begin{equation}\begin{aligned}}
\newcommand{\beas}{\begin{equation*}\begin{aligned}}
\newcommand{\eeas}{\end{aligned}\end{equation*}}
\newcommand{\eea}{\end{aligned}\end{equation}}
\newcommand{\ee}{\end{equation}}

\renewcommand{\div}{{\rm div }}

\begin{document} 
\begin{frontmatter}
\title{Strong Solutions of the Equations for Viscoelastic Fluids \\ in Some Classes of Large Data}
%%%%%
%%% algebraic decay

\author[FJ]{Fei Jiang}\ead{jiangfei0591@163.com}
% \cortext[cor1]{We  certify that the general content of this article, in whole or in part, is not submitted, accepted, or published elsewhere, including conference proceedings. }
\author[sJ]{Song Jiang}\ead{jiang@iapcm.ac.cn}
%\author[ww]{Guochun Wu}
%\ead{guochunwu@126.com}
%
\address[FJ]{College of Mathematics and
Computer Science, Fuzhou University, Fuzhou, 350108, China.}
\address[sJ]{Institute of Applied Physics and Computational Mathematics, P.O. Box 8009,  Beijing 100088, China.}
%\address[ww]{Institute of Applied Mathematics, AMSS, Chinese Academy of %Sciences, Beijing 100190, China.}

\begin{abstract}
We study the existence and uniqueness of global strong solutions to the equations of an incompressible viscoelastic fluid
in a spatially periodic domain, and show that a unique strong solution exists globally in time if the initial deformation
and velocity are small for the given physical parameters. In particular, the initial velocity can be large for the large elasticity coefficient.
The result of this paper mathematically verifies that the elasticity can prevent the formation of singularities of strong solutions
 with large initial velocity, thus playing a similar role to viscosity in preventing the formation of singularities in viscous flows.
Moreover, for given initial velocity perturbation and zero initial deformation around the rest state, we find,
 as the elasticity coefficient or time go to infinity, that
\begin{enumerate}
  \item[(1)] any straight line segment $l^0$ consisted of fluid particles in the rest state,
  after being bent by a velocity perturbation, will turn into a straight line segment that is parallel
  to $l^0$ and has the same length as $l^0$ (see Remark \ref{202003032039}).
 \item[(2)] the motion of the viscoelastic fluid can be approximated by a linear pressureless motion in Lagrangian coordinates,
  even when the initial velocity is large.
\end{enumerate}
Moreover, the above mentioned phenomena can also be found in the corresponding compressible fluid case.
   \end{abstract}
\begin{keyword}
Incompressible/compressible viscoelastic flows; strong solutions; exponential stability; elasticity coefficient.
%\MSC[2000] 35Q35\sep  76D03.
%(2000 is the default)
\end{keyword}
\end{frontmatter}

% \linenumbers

\newtheorem{thm}{Theorem}[section]
\newtheorem{lem}{Lemma}[section]
\newtheorem{pro}{Proposition}[section]
\newtheorem{cor}{Corollary}[section]
\newproof{pf}{Proof}
\newdefinition{rem}{Remark}[section]
\newtheorem{definition}{Definition}[section]
% \linenumbers

\section{Introduction}\label{introud}
\numberwithin{equation}{section}

The motion of an incompressible viscoelastic fluid can be described by the following
equations, which include a viscous stress component and a stress component for a neo-Hookean solid:
\begin{equation}\label{1.1}
\begin{cases}
\rho v_{t}+\rho v\cdot\nabla v
+\nabla p -\mu \Delta  v  = {\kappa} \mm{div} (UU^{\mm{T}})  , \\
 {U}_{t}+ v\cdot\nabla {U}= \nabla v{U} , \\
\div  v =0 ,
\end{cases}
\end{equation}
where the unknowns  ${v}:={v}(x,t)$, and ${U}:={U}(x,t)$ denote the velocity, and deformation tensor (a 3$\times$3 matrix valued function)
of the fluid, respectively. The three positive (physical) parameters $\rho$, $\mu$ and $\kappa$
stand for the density, shear viscosity coefficient and elasticity coefficient, respectively, where $\kappa$ can be defined
by the ratio between the kinetic and elastic energies, see \cite{LFHSM}, and we call the term ${\kappa} \mm{div}(UU^{\mm{T}})$ the elasticity.

For strong solutions of the both Cauchy and initial-boundary value problems for \eqref{1.1}, the authors in \cite{LFHLCZPO2,LFZPGCon2,CYZPTCPDE,LZLCZYGARMA} have established the global(-in-time) existence of solutions in various functional spaces
whenever the initial data is a small perturbation around the rest state $(0, I)$, where $I$ is the identity matrix. We also refer the reader to \cite{CPSWRRO,BJSEENARWRA,ZRFDYZTG,MNIMG,LPLMNCAMN} and references therein for the local and global existence of solutions to other closely
related models in viscoelastic fluids, and \cite{TBSMESO} for numerical evidences of singularities.
 The global existence of weak solutions to \eqref{1.1} with small perturbations near the rest state is established by Hu--Lin \cite{HXPLFHG}.
It is still, however, a longstanding open problem whether a global solution of the equations of incompressible viscoelastic fluids
 exists for any general large initial data, even in the two-dimensional case. In addition, other mathematical topics for viscoelastic fluids
 have also been widely investigated, such as the large-times behavior of solutions \cite{HXZPWHDCDSL}, the incompressible limit \cite{LZZYGE},
 the regularity of solutions \cite{CJYMNAL,KRMCTEAB} and so on.

It is well-known that viscoelasticity is a material property that exhibits both viscous and elastic characteristics
when undergoing deformation. In particular, an elastic fluid strains when stretched and quickly returns to its original state
once the stress is removed. This means that the elasticity will have a stabilizing effect in the motion of elastic/viscoelastic fluids.
Moreover, the larger $\kappa$ is, the stronger this stabilizing effect will be. Indeed, this stabilizing phenomenon
has been mathematically verified recently. For example, in \cite{LZSTCZY,LZSTCZY2,STCBT1,STCBT2} the authors used hyperbolicity of
the system \eqref{1.1} with $\mu=0$ to establish an interesting global existence result of classical solutions for small initial data
 in a subspace of $H^s(\mathbb{R}^3)$ ($s\geqslant 8$).
 %% when the initial date is also a small perturbation in that space of $(0,I)$ due to the decay property of wave spread in whole space.
 %%
 Note that such a result is not known for the Euler equations, where the elastic effect is not present. In addition,
 in \cite{JFJWGCOSdd,FJWGCZXOE,JFJSCVPDE1}, Jiang et.al. proved that the elasticity can inhibit the Rayleigh--Taylor instability
 in viscoelastic flows when $\kappa$ is properly large.

In this paper, we prove the global existence of strong solutions to the initial value problem of \eqref{1.1} defined in a periodic
domain (i.e., the motion of the viscoelastic fluid is spatially periodic), when the initial deformation (i.e, $U-I$ at $t=0$)
and the initial velocity are small for given parameters. We should remark that the initial velocity can be large if the elasticity
coefficient is large. This means that the strong elasticity can prevents the development of singularities even when the
initial velocity is large, thus playing a similar role to viscosity in preventing the formation of singularities in viscous flows.
Moreover, our analysis will reveal that for a given initial velocity perturbation and the zero initial deformation,
as $\kappa$ or $t\to\infty$,
\begin{enumerate}
  \item[(1)]
 any straight line segment $l^0$ consisted of fluid particles in the rest state, after being bent by a velocity perturbation,
  will turn into a straight line segment that is parallel to $l^0$ and has the same length as $l^0$,
 see Remark \ref{202003032039}.
  \item[(2)] roughly speaking, the motion of the viscoelastic fluid can be approximated by a linear pressureless  motion in Lagrangian
  coordinates, even though the initial velocity is large. In particular, the difference between the solutions of the linearized problem and
 those of the nonlinear problem enjoys the decay rate $\sqrt{\kappa^{-1}}e^{-\tilde{c}_3t}$, see \eqref{201905041053xxx}.
\end{enumerate}

We shall see that the above two claims will be also true for the corresponding compressible model which reads as follows.
\begin{equation}\label{1.1xx}
\begin{cases}
 \rho_t+\mm{div}(\rho{  v})=0,\\
\rho  {v}_t+\rho {v}\cdot\nabla {v}+\nabla P(\rho)= \mu\Delta v
+\lambda \nabla \mm{div}v +  \kappa\mm{div}( {UU^{\mm{T}}}/{\det U}), \\
U_t+ v\cdot\nabla U =\nabla {v} U.
\end{cases}
\end{equation}
%Next we shall complementally introduce the new mathematical notations appearing in \eqref{1.1xx}.
%
Here the pressure function $P(\tau)\in C^2(\mathbb{R})$ in \eqref{1.1xx} is always assumed to be positive and strictly increasing with respect
to $\tau$, $\lambda:=\varsigma+\mu/3$ with $\varsigma\geqslant 0$ being the bulk viscosity coefficient.
The term $\kappa \mm{div}( {UU^{\mm{T}}}/{\det U})$ is called elasticity, where the elasticity coefficient $\kappa$ is still a positive constant
for simplicity. It should be noted that ${\det U}=1$ in \eqref{1.1} due to the incompressibility, and the density $\rho$ is unknown
in \eqref{1.1xx} due to the compressibility. We mention that the well-posedness of \eqref{1.1xx} has also been widely studied, see \cite{HXPWGC,HXWDHTTIVPTCVF1,HXWDHTTIVPTCVF2,HXPGETDCMFJDE} for examples.

Since it seems difficult to directly obtain the existence of strong solutions to \eqref{1.1} for some classes of large data. Next,
we shall reformulate the motion equations in Lagrangian coordinates.

\subsection{Reformulation}

We first reformulate the motion equations \eqref{1.1} in Lagrangian coordinates, the corresponding compressible case will be dealt with in Section \ref{202001111554}.

Keeping in mind that in this paper the system \eqref{1.1} is considered in a periodic domain, without loss of generality,
we take the periodic domain to be $\mathbb{T}^3$, where $\mathbb{T}:=\mathbb{R}/\mathbb{Z}$. Let the flow map $\zeta$ be the solution to
\begin{equation}
\label{201806122101}
            \begin{cases}
\partial_t \zeta(y,t)=v(\zeta(y,t),t)&\mbox{ in }\mathbb{T}^3 \times\mathbb{R}^+ ,
\\
\zeta(y,0)=\zeta^0(y)&\mbox{ in }\mathbb{T}^3.
                  \end{cases}
\end{equation}

In Lagrangian coordinates, the deformation tensor $\tilde{U}(y,t)$ is  defined by the Jacobi matrix of $\zeta(y,t)$:
\begin{equation}\nonumber
\tilde{U}(y,t):=\nabla \zeta(y,t),\quad \mbox{i.e.},\;\; \tilde{U}_{ij}:=\partial_{j}\zeta_i(y,t) .
\end{equation}
  When we study this deformation tensor in Eulerian coordinates, it is defined  by
$$ U(x,t):=\nabla\zeta(\zeta^{-1}(x,t),t). $$
Moreover, by virtue of the chain rule, it is easy to check that $U(x,t)$ automatically satisfies the deformation equation \eqref{1.1}$_2$.
This means that the deformation tensor in Lagrangian coordinates can be directly represented by $\zeta$, if the initial data $U^0$ also satisfies
\begin{align}
\label{2020010111556}
U^0:=\nabla \zeta^0(\zeta^{-1}_0,0).
\end{align}

Next, we proceed to rewrite the elasticity in Lagrangian coordinates.
For this purpose, we introduce the matrix $\mathcal{A}:=(\ml{A}_{ij})_{3\times 3}$ and some differential operators involving $\mathcal{A}$.
Define $\ml{A}^{\mm{T}}:=(\nabla \zeta)^{-1}:=(\partial_j \zeta_i)^{-1}_{3\times 3}$. The differential operators $\nabla_{\ml{A}}$,
$\mm{div}_\ml{A}$ and $\Delta_{\ml{A}}$ are defined by $\nabla_{\ml{A}}f:=(\ml{A}_{1k}\partial_k f,
\ml{A}_{2k}\partial_kf,\ml{A}_{3k}\partial_kf)^{\mm{T}}$, $\mm{div}_{\ml{A}}(X_1,X_2,X_3)^{\mm{T}}:=\ml{A}_{lk}\partial_k X_l$
and $\Delta_{\mathcal{A}}f:=\mm{div}_{\ml{A}}\nabla_{\ml{A}}f$
for a scalar function $f$ and a vector function $X:=(X_1,X_2,X_3)^{\mm{T}}$. It should be noted that  we have used the Einstein convention of summation
over repeated indices, and $\partial_k\equiv\partial_{y_k}$.

Let $J=\det \nabla \zeta$. Obviously, $\det U |_{x=\zeta}=J$.
It is well-known that
\begin{equation}
\partial_k(J\mathcal{A}_{ik})=0.
\label{201909261909}
\end{equation}
By \eqref{201909261909} and the relation $\mathcal{A}^{\mm{T}} \nabla\zeta =I$, we have
\begin{align}
\mm{div}( {UU^{\mm{T}}}/{\det U}) |_{x=\zeta}=&
 \mm{div}_{\mathcal{A}}( \nabla \zeta \nabla \zeta^{\mm{T}}/J)=
 \mm{div}({\mathcal{A}}^{\mm{T}}(\nabla \zeta \nabla \zeta^{\mm{T}}))/J=\Delta\eta/J .  \label{201912012002}
 \end{align}

Now, assume $\det\nabla \zeta^0=1$, then $\det \nabla \zeta=1$ due to the divergence-free condition \eqref{1.1}$_3$. For the incompressible case,
 \eqref{201912012002} thus reduces to
\begin{align}
\label{202001011429}
\mm{div}( {UU^{\mm{T}}})|_{x=\zeta}=\Delta \eta.
\end{align}
Let $\eta:=\zeta-y$ and
\begin{equation*}
( u ,q)(y,t)=( v,p)(\zeta(y,t),t) \mbox{ for } (y,t)\in \mathbb{T}^3 \times\mathbb{R}^+ .
\end{equation*}
By virtue of \eqref{1.1}$_1$, \eqref{1.1}$_3$, \eqref{201806122101}$_1$ and \eqref{202001011429}, the evolution equations for $(\eta,u,q)$ in Lagrangian
coordinates read as follows.
\begin{equation}\label{01dsaf16asdfasf}
                              \begin{cases}
\eta_t=u , \\[1mm]
 {\rho} u_t-\mu \Delta_{\ml{A}}u+\nabla_{\ml{A}}q=
{\kappa}\Delta \eta, \\[1mm]
\div_\ml{A}u=0,
\end{cases}
\end{equation}
with initial data:
\begin{equation}\label{01dsaf16asdfasfsaf} (\eta, u)|_{t=0}=(\eta^0,  u^0)\;\;\mbox{ in }\mathbb{T}^3.
\end{equation}

\subsection{Results for the incompressible case}

Before stating our results, we introduce some simplified notations:
\begin{align}
& f_{\mm{h}}:=(f_1,f_2),\;\; \mathbb{R}^+_0:=[0,\infty),\;\; I_T:=(0,T),\;\; Q_T:=\mathbb{T}^3\times I_T,\;\; \int:=\int_{\mathbb{T}^3}, \nonumber \\
&\ ( w)_{\mathbb{T}^3}:=\int w\mm{d}y, \;\;  L^r:=L^r (\mathbb{T}^3)=W^{0,r}(\mathbb{T}^3),\;\; {H}^k:=W^{k,2}(\mathbb{T}^3),\nonumber \\
& \|\cdot \|_{k} :=\|\cdot \|_{H^k},\;\; I_0^{\mm{h}}( w,\eta):=\|(w,\sqrt{\kappa} \nabla \eta) \|_2^2,\nonumber \\
& H^3_{*}:=\{\eta\in H^3~|~ \eta(y, t)+y: \mathbb{R}^3 \to \mathbb{R}^3 \mbox{ is a } C^1\mbox{-diffeomorphism mapping}\}\nonumber  ,\\
& H^3_{*,1}:=\{\eta\in H^3_*~|~\det\nabla(\eta +y )=1\} ,\;\; H^1_\sigma:=\{u\in H^1~|~\mm{div}u=0\},  \nonumber \\ &
 \underline{H}^k:=\{w\in H^k~|~( w)_{\mathbb{T}^3}=0\} , \;\; H^k_\sigma:=H^1_\sigma \cap H^k,\nonumber \\
&  \mathfrak{U}_T:=\{v\in C^0([0,T), H^2)~|~ v_t\in L^\infty(I_T, L^2), \;\;  \nabla v\in L^2(I_T,H^2),\;\; v_t\in L^2(I_T,H^1)\}, \nonumber \\
 & \|v\|_{\mathfrak{U}_T}:=\sqrt{\|v\|_{C^0(\overline{I_T}, H^2)}^2 +\|v_t  \|_{L^\infty(I_T, L^2)}^2+ \|( v,v_t )\|_{L^2(I_T,H^3\times H^1)}^2} ,\nonumber\\
& \|{q}^k\|_{\mathfrak{Q}_T}:=  \sqrt{\|\nabla q\|_{L^\infty(I_T, L^2 )}^2 +  \|\nabla q\|_{L^2(I_T, H^1 )}^2}  ,\nonumber
 \end{align}
where $1< r\leqslant \infty$ and  $k$ is a non-negative integer.

The letters ${c}_i$ respectively $ {c}_i^{\mm{P}}$ ($1\leqslant i \leqslant 5$) are fixed constants which may depend on
the parameters $\rho$ and $\mu$ respectively $P(\cdot)$, $\bar{\rho}$, $\mu$ and $\lambda$ in \eqref{n0101nnn}$_2$. The letters
$c_0$, $c$, $c_{\mm{P}}$ and $c_\kappa$ will denote generic positive constants that may vary from line to line. Moreover, $c_0$ is independent of any parameters, while
$c_{\kappa}$ may depend on $\rho$ and $\mu$ as well as $\kappa$, $c$ and $c_{\mm{P}}$ on some parameters such as ${c}_i$ and $ {c}_i^{\mm{P}}$ respectively.
$a\lesssim_0 b$, $a\lesssim b$, $a\lesssim_{\mm{P}} b$ and $a\leqslant c_\kappa b$ mean that $a\leqslant c_0b$,
$a\leqslant cb$, $a\leqslant c_{\mm{P}} b$ and $a\lesssim_{\kappa} b$, respectively.

Next, we state the main results for the incompressible case. The first result is concerned with the existence of strong solutions to the initial value problem
\eqref{01dsaf16asdfasf}--\eqref{01dsaf16asdfasfsaf} in some classes of large initial data:
%%%%%%%%%%%%%%%%%%%%%%%%%%%%%%%%%%%%%%%%
\begin{thm}
\label{201904301948}
 There are constants ${c}_1\geqslant 1$ and ${c}_2\in (0,1]$, such that for any $(\eta^0,u^0)\in H^{3}_{*,1}\times H^{2}$ and $\kappa$, satisfying
$\mm{div}_{\mathcal{A}^0}u^0=0$ and
\begin{equation}
\label{201909281832} \kappa\geqslant \frac{1}{c_2}\max\left\{ 2 \sqrt{ {c}_1I_0^{\mm{h}}(u^0,\eta^0)} ,( 4 {c}_1I_0^{\mm{h}}(u^0,\eta^0))^2 \right\},
\end{equation}
where ${\mathcal{A}^0}$ denotes the initial data of $\mathcal{A}$,
the initial value problem \eqref{01dsaf16asdfasf}--\eqref{01dsaf16asdfasfsaf} admits a unique global solution
$(\eta, u,q)\in C^0(\mathbb{R}_0^+$, $H^{3})\times\mathfrak{U}_\infty\times L^\infty(\mathbb{R}^+, \underline{H}^2 )$.
Moreover, the solution $(\eta,u)$ enjoys the stability estimate
\begin{align}
 \|( u, \sqrt{\kappa}\nabla \eta)\|_2^2 +\int_0^t \| \nabla (u,\sqrt{\kappa}\eta)\|_2^2 \mm{d}\tau \lesssim I_0^{\mm{h}}(u^0,\eta^0).\label{201905041053}
\end{align}
In addition, it holds that
\begin{align}& \|\div \eta \|^2_2\lesssim \|\nabla \eta\|_2^2 ,\label{202001121649} \\
& \| \nabla q\|^2_1 \lesssim\|\nabla u\|_1^2+{\kappa}\|\nabla \eta\|_2^2 , \label{202001071328} \\
&\|\tilde{\mathcal{A}}\|_2 \lesssim_0 \|\nabla \eta\|_2 \lesssim_0 1,  \label{202001121644}
 \end{align}
 where $\tilde{\mathcal{A}}:=\mathcal{A}-I$.
\end{thm}

Now we briefly sketch the proof idea of Theorem \ref{201904301948}. If we multiply \eqref{01dsaf16asdfasf}$_2$ by $u$ in $L^2$, then we get the basic energy identity:
\begin{align}
\label{ssebdaiseqinM0846dfgssgsdxx}
\frac{1}{2}\frac{\mm{d}}{\mm{d}t}\left(\rho\|u\|_0^2+ {\kappa} \|\nabla \eta\|_0^2\right)
+\mu \|\nabla_{\mathcal{A}} u\|_0^2= 0,
\end{align}
which implies
\begin{align}
\label{ssebdsdxx}
 \frac{\rho}{ \kappa}\|u\|_0^2+ \|\nabla \eta\|_0^2
+\frac{2\mu}{\kappa}\int_0^t \|\nabla_{\mathcal{A}} u\|_0^2\mm{d}\tau\leqslant
\frac{I_0}{\kappa},
\end{align}
where $I_0:={\rho} \|u^0\|_0^2 +\kappa\|\nabla \eta^0\|_0^2$.
In particular, we see  that $\|\nabla \eta\|_0\to 0$ as $\kappa\to \infty$ for fixed $I^0$. This key observation motivates us to conclude that the equation \eqref{01dsaf16asdfasf}
may be approximated by the following linear equations for sufficiently large $\kappa$:
\begin{equation}\label{01dsaf16asdfasfxx}
                              \begin{cases}
\eta_t^{\mm{l}}=u^{\mm{l}}, \\[1mm]
 {\rho} u_t^{\mm{l}}-\mu \Delta u^{\mm{l}}+\nabla q^{\mm{l}}=
{\kappa}\Delta \eta^{\mm{l}}, \\[1mm]
\div u^{\mm{l}}=0.
\end{cases}
\end{equation}
Since the linear equations have global solutions with large initial data, we could expect that the initial value problem \eqref{01dsaf16asdfasf}--\eqref{01dsaf16asdfasfsaf}
may also admit a global large solution for sufficiently large $\kappa$ (with fixed $I_0^{\mm{h}}(u^0,\eta^0)$) as stated in Theorem \ref{201904301948}.

Thus, to obtain Theorem \ref{201904301948}, the key step is to derive the \emph{a priori} estimate \eqref{201905041053}
under sufficiently large $\kappa$. Fortunately, by careful energy estimates, we find that the estimate \eqref{ssebdsdxx}
still holds for higher order spatial derivatives of $(\eta,u)$ under sufficiently large $\kappa$. Namely,
we can conclude that there are constants $K$ (depending possibly on $I_0^{\mm{h}}(u^0,\eta^0)$ but not on $T$) and $\delta$, such that
\begin{equation}\label{aprpiose1}
\sup_{ 0\leqslant  t <T  } \| (u(t), \sqrt{\kappa}\nabla \eta(t))\|_2 \leqslant  {K}/ {2},
\end{equation}
provided that
\begin{equation}\label{aprpiosesnewxxxx}
 \sup_{  0\leqslant  t<T} \|(u(t),\sqrt{\kappa}\nabla \eta (t))\|_2\leqslant K \;\; \mbox{ for any given }\; T
\end{equation}
and
\begin{equation}\label{aprpiose1snewxxxxz}
\max\{K, K^4\}/\kappa\in (0,\delta^2].\end{equation}

Based on the above fact and the existence of a unique local solution, we can immediately obtain Theorem \ref{201904301948}.
The detailed proof will be presented in Section \ref{201805081dfs731}. In addition, the proof of the existence of a unique local solution
will be given in Section \ref{201912062021}.

The second result is concerned with the properties of the solution $(\eta,u)$ to the initial value problem \eqref{01dsaf16asdfasf}--\eqref{01dsaf16asdfasfsaf} in Theorem \ref{201904301948}.
\begin{thm}\label{201912041028}
Let $(\eta,u)$ be the solution established in Theorem \ref{201904301948}, then we have
 \begin{enumerate}
   \item[(1)] Exponential stability of $(\eta,u)$:
  \begin{align}
& \|(\bar{u},\sqrt{\kappa}\nabla \eta)\|_2^2 +\int_0^t( \| \bar{u}\|_3^2+{\kappa}\|\nabla \eta\|_2^2 )e^{ {c}_3(\tau-t)}\mm{d}\tau
\lesssim e^{- {c}_3 t}I_0^{\mm{h}}(\bar{u}^0,\eta^0),\quad \forall\, t\geqslant 0,\label{201905041x053}
\end{align}
where $\bar{u}:=u-({u}^0)_{\mathbb{T}^3}$  and $\bar{u}^0:=u^0-({u}^0)_{\mathbb{T}^3}$.
 \item[(2)] Large-time behavior of $\eta$:
    \begin{align}
 & \|\bar{\eta}\|_{2}^2 + \int_0^t \|\bar{\eta}\|_2^2\mm{d}\tau
  \lesssim e^{-  {c}_3 t}I_0^{\mm{h}}(\bar{u}^0,\eta^0)/{\kappa},\label{201905041053x}\\
&  \| \eta-( {u}^0)_{\mathbb{T}^3}t-\varpi\|_3^2 +\int_0^t  \| \eta-( {u}^0)_{\mathbb{T}^3}t-\varpi\|_3^2 e^{ {c}_3(\tau-t)}\mm{d}\tau \lesssim e^{- {c}_3 t}I_0^{\mm{h}}(\bar{u}^0,\eta^0)/{\kappa},\label{20191211197}
  \end{align}
where $\bar{\eta}:=\eta-\eta(y_{\mm{h}},0 ,t)$ and $\varpi=(\eta^0)_{\mathbb{T}^3}/8\pi^3$.
   \item[(3)] Stability of $(\eta,u)$ around $( \eta^{\mm{l}}, u^{\mm{l}})$:
 \begin{align}
& \| u^{\mm{d}}\|_2^2+ {\kappa}\| \eta^{\mm{d}}\|_3^2+  \int_0^t \| (u^{\mm{d}}, \sqrt{\kappa}\eta^{\mm{d}})\|_3^2 e^{ {c}_3(\tau-t)}
\mm{d}\tau\nonumber\\
 &\quad \lesssim \sqrt{\kappa^{-1}}\max\{1,\sqrt{\kappa}^{-1}\}e^{- {c}_3 t}\nonumber \\
&\qquad \times \left(\sqrt{I_0^{\mm{h}}(u^0,\eta^0)} +I_0^{\mm{h}}(u^0,\eta^0) \right)
 {( \| \nabla \eta^0 \|_2^2 I_0^{\mm{h}}(  {u}^0, {\eta}^0) +I_0^{\mm{h}} ( \bar{u}^0, {\eta}^0))}.
\label{201905041053xxx}
\end{align}
  Here $(\eta^{\mm{d}}, u^{\mm{d}} ):= (\eta-\eta^{\mm{l}} , u-u^{\mm{l}} )$ and
  $(\eta^{\mm{l}}, u^{\mm{l}} )\in C^0(\mathbb{R}_0^+$, $H^{3})\times \mathfrak{U}_\infty $ is the unique strong solution
  of the linear system \eqref{01dsaf16asdfasfxx} with $q^{\mm{l}}=0$ and initial data:
\begin{equation}
(\eta^{\mm{l}},u^{\mm{l}})|_{t=0}=(\eta^0+\eta^{\mm{r}}, u^0+u^{\mm{r}}),
\label{202001070914}
\end{equation}
where $(\eta^{\mm{r}}, u^{\mm{r}})\in  \underline{H}^3\times \underline{H}^2$ has the following properties:
\begin{enumerate}[(a)]
  \item $\mm{div}(u^0+u^{\mm{r}})=0$ and $\|u^{\mm{r}} \|_2\lesssim_0 \|\nabla \eta^0\|_2\|u^0\|_2$.
  \item  $\mm{div}(\eta^0+\eta^{\mm{r}})=0$  and $\|\eta^{\mm{r}} \|_3\lesssim_0 \|\nabla \eta^0\|_2^2$.
\end{enumerate}
 \end{enumerate}
 \end{thm}
\begin{rem}
We give an explanation why the initial data $(\eta^0,u^0)$  has to be modified as in \eqref{202001070914}.
\begin{enumerate}
  \item[(1)]
Since the initial data for $u^{\mm{l}}$ has to satisfy the divergence-free condition, i.e., $\mm{div}(u^{\mm{l}}|_{t=0})=0$,
one thus has to adjust the initial data $u^0$ as in \eqref{202001070914}.
  \item[(2)]
 The initial data $\eta^0$ for $\eta$ can be directly used as an initial data for $\eta^1$. In this case one can see that $\mm{div}\eta^{\mm{l}}=\mm{div}\eta^0$. Consequently, the decay-in-time of $\nabla \eta^{\mm{d}}$ by \eqref{201905041x053}
  can not be expected, unless $\mm{div}\eta^0=0$. Hence, we have to modify $\eta^0$ as in \eqref{202001070914}, so that
  the new initial data``$\eta^0+\eta^{\mm{r}}$'' also satisfies the divergence-free condition.
\end{enumerate}
\end{rem}
%%%%%%%%%%%%%%%%%%%%%%%%%%%%%%%%%%%%%%%
\begin{rem}
\label{202001062012} Let us try to give a physical meaning hidden in \eqref{201905041053x}. Consider the viscoelastic fluid
in a periodic cell $(0,2\pi)^3$, and think that the fluid is made up of infinite fluid segments that are parallel
to $x_3$-axis. Now, we consider a  straight line segment denoted by $l^0$ in $(0,2\pi)^3$ at $t=0$, and any given (fluid) particle
$y$ in $l^0$. The two particles at the upper and lower endpoints are denoted by $y^1$ and $y^2$, and note that
$y_{\mm{h}}=y_{\mm{h}}^1=y_{\mm{h}}^2$. We disturb the rest state by a perturbation $(\eta^0,u^0)$ at $t=0$. Then, the line segment $l^0$ will be bent and move to a new location at time $t$ that may be curved line segment and denoted by $l^t$.
% (it is a curve when $l^0$ is bent at time $t$).
Since the motion of the fluid is spatially periodic,  the segment $y^1y^2$ is thus parallel to $x_3$-axis
 and can be represented by
$$l^{\mm{n}}: x_{\mm{h}}=\eta_{\mm{h}}(y_{\mm{h}},0,t)+y_{\mm{h}}, \ x_3=\eta_3(y_{\mm{h}},0,t)+y_3,\ 0\leqslant y_3\leqslant 2\pi .$$
At time $t$, the new location of a particle $y$ on $l^0$ is given by $\eta(y,t)+y$.
Thus we see that $|\bar{\eta} (y,t)|$ has a geometric meaning, namely, it represents the distance between the two points
$\eta(y,t)+y$ on $l^t$ and $(\eta_{\mm{h}}(y_{\mm{h}},0,t)+y_{\mm{h}}, \eta_3(y_{\mm{h}},0,t)+y_3)$ on $l^{\mm{n}}$.

By \eqref{201905041053x} and the interpolation inequality \eqref{201905041149}, we see that
$$|\bar{\eta} (y,t)|\lesssim e^{-{c}_3 t} I_0^{\mm{h}}(\bar{u}^0,\eta^0) /\kappa.$$
In particular, for the case $\eta^0=0$,
 $$|\bar{\eta} (y,t)|\lesssim  e^{-  {c}_3 t}   \|u^0\|_2^2 /\kappa   ,$$
from which and the geometric meaning of $|\bar{\eta} (y,t)|$ we immediately see that the curve $l^t$
oscillates around $l^{\mm{n}}$, and the amplitude tends to zero, as $\kappa$ or $t$ is extremely large.
\end{rem}
\begin{rem}
\label{202003032039}
Similar to Remark \ref{202001062012}, we can also give a physical meaning hidden in \eqref{20191211197}. Namely,
consider a straight line segment $l^0$ consisted of fluid particles in the rest state and perturb
the segment by a velocity, thus it will be bent. Then, as $\kappa$ or $t\to\infty$,
the perturbed segment will turn into a straight line segment that is parallet to $l^0$ and has the same length as $l^0$.
\end{rem}

The detailed derivation of Theorem \ref{201912041028} will be given in Section \ref{201912042253}.

Finally, we mention that the existence result in Lagrangian coordinates from Theorem \ref{201904301948} can be recovered
to the one in Eulerian coordinates. In fact, noting that the solution $\eta$ in Theorem \ref{201904301948} satisfies
\begin{align}&
   \zeta:= \eta(y, t)+y   : \mathbb{R}^3 \to \mathbb{R}^3 \mbox{ is a }C^1\mbox{ diffeomorphism mapping}, \label{2312018031adsadfa21601xx}
\end{align}
and using an inverse transform of the Lagrangian coordinates, i.e., $(v,U,p):=(u,\nabla \zeta,q)|_{y=\zeta^{-1}}$, we can easily
get a global solution of \eqref{1.1} from Theorem \ref{201904301948}. More precisely, we have the following conclusion,  the proof of which
is referred to \cite[Theorem 1.2]{JFJSJMFMOSERT}.
%%%%%%%%%%%%%%%%%%%%%%%%%%%%%%%%%%%%%%%%%%%%%%%%%%%%%%%%%%%%%%%%%%%%%%%%%%%%
\begin{thm}\label{201904301948xx}
Let  $ (v^0,U^0)\in H^2_\sigma \times H^{2} $  and $\kappa$ satisfy
 \begin{enumerate}[(1)]
  \item  $U^{0}=\nabla\zeta^{0}(\zeta_{0}^{-1}(x))$, where $\eta^{0}:=\zeta^{0}(y)-y\in H^3_{*,1}$,
\item $\mm{div}_{\mathcal{A}^0}(v^0(\zeta^0))=0$,
  \item \eqref{201909281832} holds with $v^0(\zeta^0)$ in place of $u^0$.
\end{enumerate}
Then, the initial value problem of \eqref{1.1} with initial data $(v,U)|_{t=0}=(v^0,U^0)$ admits a unique global solution
$((v,U),p)\in C^0(\mathbb{R}_0^+,H^2_\sigma\times H^2)\times  L^\infty(\mathbb{R}^+, \underline{H}^2 )$.
 \end{thm}

\subsection{Results for the compressible case}\label{202001111554}

In this subsection we describe a further extension of the above results for incompressible viscoelastic fluids to the corresponding compressible fluid case.

Similarly to \eqref{01dsaf16asdfasf}, we first rewrite \eqref{1.1xx} in Lagrangian coordinates. For this purpose, let $\zeta$ be the flow map as in \eqref{201806122101}.
If $U^0$ satisfies \eqref{2020010111556}, similarly to \eqref{01dsaf16asdfasf}, we can easily write the evolution equations for $(\eta,\varrho,u)$ in Lagrangian coordinates
as follows.
\begin{equation}\label{n0101nn}
\begin{cases}
\eta_t=u  ,\\[1mm]
\varrho_t+\varrho\mm{div}_{\ml{A}}u=0 ,\\[1mm]
{\varrho} u_t+ \nabla_{\ml{A}} P(\varrho)= \mu  \Delta_{\ml{A}} v
+\lambda \nabla_{\ml{A}} \mm{div}_{\ml{A}}v +\kappa \Delta \eta /J ,\\[1mm]
\end{cases}\end{equation}
where $\varrho (y,t):= \rho (\eta(y,t)+y,t)$.

 Then, we consider a rest state $(\bar{\rho},0,I)$ of \eqref{1.1xx} and expect the time-asymptotic behavior like
\begin{equation}
\label{202001062252}
\zeta(y,t)\to y\mbox{ and } \varrho(y,t)\to \bar{\rho}\;\; \mbox{ as }t\to \infty.
\end{equation}
It follows from \eqref{n0101nn}$_1$ that
\begin{equation}\label{Jtdimau}
J_t=J\mm{div}_{\mathcal{A}}u,             \end{equation}
which, together with \eqref{n0101nn}$_2$, yields
\begin{equation}\label{0122sd} \partial_t(\varrho J)=0.   \end{equation}
Thus we deduce from the asymptotic behavior \eqref{202001062252} and \eqref{0122sd} that $\varrho_0J_0=\varrho J=\bar{\rho}$,
which implies $\varrho = \bar{\rho}J^{-1}$, provided the initial data $(\varrho_0,J_0)$ satisfies
\begin{align}
&\label{abjlj0i}  \varrho_0 =\bar{\rho}J^{-1}_0 ,
\end{align}
where $\varrho_0$ and $J_0$ are the initial data of $\varrho$ and $J$, respectively.

Therefore, if one further imposes the assumption \eqref{abjlj0i} initially, the equations \eqref{n0101nn} can be written as follows.
\begin{equation}\label{n0101nnn} \begin{cases}
\eta_t=u ,\\[1mm]
\bar{\rho}  u_t+ J(\nabla_{\ml{A}} P(\bar{\rho}J^{-1})- \mu \Delta_{\ml{A}}u
-\lambda \nabla_{\ml{A}} \mm{div}_{\ml{A}}u)=\kappa \Delta \eta.
\end{cases} \end{equation}
Finally, we impose the initial value conditions for \eqref{n0101nnn}:
\begin{equation}  \label{2201912032006}
 (\eta ,u)|_{t=0}=(\eta^0,u^0).
 \end{equation}
Then we have the following conclusion for the initial value problem \eqref{n0101nnn} and \eqref{2201912032006}:
%%%%%%%%%%%%%%%%%%%%%%%%%%%%%%%%%%%%%%%%%%%%%%%%%%%%
\begin{thm}   \label{x2019saf04301948}
 There are constants ${c}_1^{\mm{P}}\geqslant 1$ and $ {c}_2^{\mm{P}}\in (0,1]$, such that for any $(\eta^0,u^0)\in H^{3}_{*}\times H^{2}$ and $\kappa$ satisfying
\begin{equation} \nonumber
 \kappa\geqslant \frac{1}{ {c}_2^{\mm{P}}}\max\left\{ \left(2\sqrt{ {c}_1^{\mm{P}} I_0^{\mm{h}}(u^0,\eta^0) } \right)^{2/3},( 4 {c}_1^{\mm{P}} I_0^{\mm{h}}(u^0,\eta^0))^2 \right\},
\end{equation}
 the initial  value problem \eqref{n0101nnn} and \eqref{2201912032006} admits a unique global solution $(\eta, u)\in C^0(\mathbb{R}_0^+,H^{3}_*\times H^2)$.
Moreover, the solution $(\eta,u)$ enjoys the stability estimate \eqref{201905041053} with ``$\lesssim_{\mm{P}}$" in place of ``$\lesssim $", \eqref{202001121644} and  $1/ 2 \leqslant \det (\nabla \eta+I)\leqslant 3/2$.
\end{thm}

The properties on the solution given in Theorem \ref{201912041028} can also be extended to the solution $(\eta,u)$ established in Theorem \ref{x2019saf04301948}. Namely,
\begin{thm}\label{20191asfas2041028}
Let $(\eta,u)$ be the solution given by Theorem \ref{x2019saf04301948}, then
 \begin{enumerate}
   \item[(1)]  $(\eta,u)$ satisfies \eqref{201905041x053}--\eqref{20191211197} with some constant $ {c}_3^{\mm{P}}$ and ``$\lesssim_{\mm{P}}$''
        in place of $ {c}_3$ and  ``$\lesssim $'', respectively.
   \item[(3)] The following stability of $(\eta,u)$ around the solution $( \eta^{\mm{l}},u^{\mm{l}})$ holds.
 \begin{align}
 & \| u^{\mm{d}}\|_2^2 + {\kappa}\| \eta^{\mm{d}}\|_3^2+\int_0^t\| (u^{\mm{d}},\sqrt{\kappa}\eta^{\mm{d}})\|_3^2
 e^{ {c}_3^{\mm{P}}(\tau-t)}\mm{d}\tau  \nonumber   \\
 & \qquad \qquad \lesssim_{\mm{P}} \sqrt{\kappa}^{-1}\max\{1,\kappa^{-1}\} e^{-{c}_3^{\mm{P}} t}|I_0^{\mm{h}}(u^0,\eta^0)|^{3/2},
 \label{201905safa041053xxx}
\end{align}
 where $(\eta^{\mm{d}}, u^{\mm{d}}):= (\eta-\eta^{\mm{l}} , u-u^{\mm{l}})$, $(\eta^{\mm{l}},u^{\mm{l}})\in C^0(\mathbb{R}_0^+,H^{3})\times \mathfrak{U}_\infty $
 is the unique strong solution of the following linearized problem:
 \begin{equation}\label{01dsaf16asdfsadffasfxx}
  \begin{cases}
\eta_t=u, \\[1mm]
\bar{\rho} u_t-\nabla( P'(\bar{\rho})\bar{\rho}\,\mm{div}\eta+\lambda\mm{div}u )-\Delta(\mu  u +{\kappa}\eta)=0,\\
 (\eta,   u)|_{t=0}=(\eta^0,  u^0).
\end{cases}
\end{equation}
 \end{enumerate}
 \end{thm}
 %%%%%%%%%%%%%%%%%%%%%%%%%%%%%%%%%%%%%%%%%%%%
\begin{rem}
Analogously to Theorem \ref{201904301948xx}, one can also obtain an existence result of a unique global solution to \eqref{1.1xx} from Theorem \ref{x2019saf04301948}.
\end{rem}

The proof of Theorems \ref{x2019saf04301948} and \ref{20191asfas2041028} is similar to that for the incompressible case and will be sketched in Section \ref{20201011926}.

\section{Proof of Theorem \ref{201904301948}}\label{201805081dfs731}

This section is devoted to the proof of Theorem \ref{201904301948}, in which the key step is to establish the (\emph{a priori}) stability estimate \eqref{201905041053}
under a proper assumption. To this end, let $(\eta,u,q)$ be a solution of the initial value problem \eqref{01dsaf16asdfasf}--\eqref{01dsaf16asdfasfsaf} defined on $Q_T$ for any given $T>0$,
where $\eta^0 \in H^{3}_{*,1}$, $\mm{div}_{\mathcal{A}^0}u^0=0$ and $(q)_{\mathbb{T}^3}=0$.

We further assume that $(\eta,u)$ and $K$ satisfy \eqref{aprpiosesnewxxxx} and \eqref{aprpiose1snewxxxxz},
 where $K\geqslant 1$ will be defined by \eqref{201911262060}, and $\delta\in (0,1]$ is a sufficiently small constant.
 We should keep in mind that the smallness of $\delta$ only depends on the parameters $\rho$ and $\mu$
%%%% and will be repeatedly used
 in what follows. In addition, we see from \eqref{aprpiosesnewxxxx} and \eqref{aprpiose1snewxxxxz} that
\begin{equation}\label{aprpiosesnewxxxx1x}
 \sup_{0\leqslant t < T} \|\nabla \eta(t)\|_2 \lesssim_0\delta .
\end{equation}

Next our aim is to establish the stability estimate \eqref{201905041053}, while $(\eta,u)$ enjoys the estimate \eqref{aprpiose1}.
To begin with, we introduce some preliminary estimates.
\begin{lem}    \label{201805141072}
\begin{enumerate}
  \item[(1)] Let $\eta$ satisfy \eqref{aprpiosesnewxxxx1x} with $\delta\in (0,1]$, then
\begin{align}
&\label{aimdse}
 \|\tilde{\mathcal{A}}\|_2 \lesssim_0 \|\nabla \eta\|_2 \lesssim_0 \delta,\\
& \|\tilde{\mathcal{A}}_t\|_1\lesssim_0 \| \nabla u\|_1 \label{201905031644},\\
 & \label{2016957x} \|\div \eta \|_2\lesssim_0 \|\nabla \eta\|_2^2.
\end{align}
  \item[(2)] Let $\eta$ satisfy \eqref{aprpiosesnewxxxx1x} with $\delta$ sufficiently small, then
 \begin{align}
& \label{20160614fdsa1957}
\|\nabla w\|_0 \lesssim_0 \|\nabla_{\mathcal{A}}w\|_0\lesssim_0 \|\nabla w\|_0\quad\mbox{for any }\nabla w\in L^2.
\end{align}
\item[(3)] It holds that for any given $i\geqslant 1$,
\begin{align}& \label{20160614fdsa1957x}
\|\nabla^i w\|_0  \leqslant  \| w\|_i \leqslant \tilde{c}_i \|\nabla^i w\|_0\quad\mbox{for any }  w\in  \underline{H}^i,
\end{align}
\end{enumerate}
where the constant $\tilde{c}_i>0$ only depends on $i$.
\end{lem}
\begin{pf}
Let us recall some well-known inequalities.
\begin{itemize}
\item Embedding inequality: for given $1\leqslant r\leqslant 6$,
\begin{equation}
\label{201905041124}
\|f\|_{L^r}\leqslant \tilde{c}_r \|f\|_1\;\;\mbox{ for any }f\in H^1,
\end{equation}
where the constant $\tilde{c}_r$ only depends on $r$.
  \item   Interpolation  inequalities: let $1\leqslant j<i$, then
\begin{alignat}{2}
& \label{201905041149}
\|f\|_{L^\infty}\lesssim_0  \|f\|_{1}^{1/2}\|f\|_2^{1/2}& & \;\; \mbox{ for any }f\in H^2,\\
 &\label{201905041130}
 \| f\|_j\leqslant \tilde{c}_{i,j}  \|  f\|_0^{1-j/i}\|f\|_i^{j/i}& &\;\; \mbox{ for any }f\in H^i,
\end{alignat}
where  the  constant $\tilde{c}_{i,j}$ only depends on $i$ and $j$.
  \item
  Poinc\'are's inequality: let $r\geqslant 1$, then
  $$\left\|f\right\|_{L^r}^r\leqslant\tilde{c}_r\left(\|\nabla f\|_{L^r}^r +\left|\int f\mm{d}y\right|^r\right)\;\;\mbox{ for any } f\in W^{1,r},$$
where the constant $\tilde{c}_r$ only depends on $r$.
\end{itemize}
%% Next we will derive the desired estimates in Lemma \ref{201805141072}.

(1) Noting that $\det(\nabla\eta^0 +I)=1$, then
\begin{equation}   \label{201812081507}
\det(\nabla \eta+ I)=1.
\end{equation}
Recalling the definition of $\mathcal{A}$, we see that
\begin{equation}   \label{06131422}
\mathcal{A}=(A^*_{ij})_{3\times 3},
\end{equation}
where $A^{*}_{ij}$ is the algebraic complement minor of $(i,j)$-th entry of the matrix $(\partial_j \zeta_i)_{3\times 3}$.

 Using H\"older's inequality, one can easily derive from \eqref{06131422} that
\begin{align}
 \|\tilde{\mathcal{A}}\|_2 \lesssim_0 (1+\|\nabla \eta\|_{L^\infty}) \| \nabla \eta\|_2 +
  \|\nabla^2\eta \|_{L^3}\| \nabla^2 \eta\|_{L^6}. \label{202001192046}
   \end{align}
If we make use of \eqref{aprpiosesnewxxxx1x}, \eqref{201905041124} and \eqref{201905041149}, we immediately get \eqref{aimdse} from \eqref{202001192046}.
In a similar manner, we can also derive \eqref{201905031644} by employing \eqref{01dsaf16asdfasf}$_1$.

   In view of \eqref{201812081507}, we easily verify that
\begin{align}
\label{202001081820}
1=\det (\nabla \eta+ I)=1+ \mm{div} \eta+r_\eta,
\end{align}
where $r_\eta:= ((\mm{div}\eta)^2-\mm{tr} (\nabla \eta)^2)/2+ \det \nabla \eta$. Thus, $\mm{div} \eta=-r_\eta$, which implies that
\begin{align}
 \|\mm{div}\eta\|_2 \lesssim_0 & (1+\|\nabla \eta\|_{L^\infty})(\|\nabla \eta\|_{L^\infty} +   \| \nabla \eta \|_{2} ) \|\nabla \eta\|_2  \lesssim_0  \| \nabla \eta\|_2^2.
 \label{2017020614181721xxdfasdfasadfaff}
   \end{align}
Therefore, \eqref{2016957x} is obtained.

(2) Exploiting \eqref{aimdse} and \eqref{201905041149}, one has
\begin{align}  &
 \|\nabla_{\mathcal{A}}w\|_0\lesssim_0 (1+\|\tilde{\mathcal{A}}\|_{L^\infty}) \|\nabla w\|_0 \lesssim_0 \|\nabla w\|_0 , \nonumber  \\
 &   \| \nabla w\|_0\lesssim_0 \|\nabla_{\tilde{\mathcal{A}}}w\|_0+ \|\nabla_{{\mathcal{A}}} w\|_0\lesssim_0  \delta\|\nabla w\|_0 + \|\nabla_{\mathcal{A}}w\|_0,\nonumber
\end{align}
which imply \eqref{20160614fdsa1957} for sufficiently small $\delta$.

(3) Finally, \eqref{20160614fdsa1957x} easily follows from Poinc\'are's inequality. This completes the proof of Lemma \ref{201805141072}.
 \hfill $\Box$
\end{pf}
%%%%%%%%%%%%%%%%%%%%%%%%%%%%%%%%%%%%%%%%%%%%%%%%%%%%%%
\begin{lem}\label{201612132242nx}
Under the condition \eqref{aprpiosesnewxxxx1x} with sufficiently small $\delta$, the following estimate hold.
\begin{align}
  & \|\nabla q\|_1 \lesssim   \|\nabla u\|_1^2 + {\kappa}   \|\nabla \eta\|_2^2.  \label{2017020614181721}
   \end{align}
\end{lem}
\begin{pf} Let $\alpha$ satisfy $|\alpha|=1$. Applying $\partial^\alpha\mm{div}_{\mathcal{A}}$ to \eqref{01dsaf16asdfasf}$_2$,
and using the identity $\mm{div}_{\mathcal{A}} u_t = -\mm{div}_{\mathcal{A}_t} u $, we have $ \partial^\alpha\Delta  q = \partial^\alpha f$, where
$$f:={\kappa}(\Delta\mm{div}\eta+ \mm{div}_{\tilde{\mathcal{A}}}\Delta \eta) +\rho \mm{div}_{\mathcal{A}_t}u
- (\mm{div}_{\tilde{\ml{A}}}\nabla_{\ml{A}}q + \mm{div} \nabla_{\tilde{\ml{A}}}q) .$$

From the regularity theory of elliptic equations one gets
\begin{align}
& \| \nabla q\|_1 \lesssim_0 \|f\|_0.  \label{201711118xx41}
\end{align}
In addition, it is easy to deduce that
\begin{align}
\|  f \|_0 \lesssim & \kappa(   \| \Delta \mm{div}\eta \|_{0}
 + \|\tilde{\mathcal{A}}\|_{L^\infty} \|\nabla^3\eta\|_0 ) + \|\mathcal{A}_t\|_{L^3}\|\nabla u\|_{L^6} \nonumber  \\
 & + \|\nabla \mathcal{A}\|_{L^3}\|\nabla q \|_{L^6}(1+\|\tilde{\mathcal{A}}\|_{L^\infty})+\| \tilde{\mathcal{A}} \|_{L^\infty}
   \|\nabla^2q\|_0(1+\| {\mathcal{A}} \|_{L^\infty})\nonumber \\
 \lesssim & \kappa(   \|  \mm{div}\eta \|_2  +   \|\nabla\eta\|_2^2 ) + \|\mathcal{A}_t\|_{1}\|\nabla u\|_{1}
 + \delta(\|\nabla  q \|_{L^6} +  \|\nabla^2q\|_0)\nonumber  \\
\lesssim & \kappa  \| \nabla \eta\|_2^2  +  \|\nabla u\|_1^2 +  \delta\|\nabla^2 q \|_{0}, \label{201909121931}
   \end{align}
where we have used \eqref{aimdse}, \eqref{201905041124} and \eqref{201905041149} in the second inequality, and
 \eqref{201905031644},  \eqref{2016957x} and Poinc\'are's inequality in the last inequality.
Finally, inserting \eqref{201909121931} into \eqref{201711118xx41}, we get \eqref{2017020614181721}. \hfill $\Box$
\end{pf}

Next we proceed to derive more estimates of $(\eta,u)$.
\begin{lem}\label{201612132242nn}
Under the assumptions \eqref{aprpiosesnewxxxx}--\eqref{aprpiose1snewxxxxz} with sufficiently small $\delta$, it holds that
\begin{align}
& \frac{\mm{d}}{\mm{d}t}  \left(  {\rho} \|\nabla^2 u\|_0^2+ \kappa \| \nabla^3\eta \|_0^2\right)+
  \mu \|\nabla^3 u \|_0^2   \lesssim
\delta \kappa \|\nabla \eta\|_2^2.
\label{ssebdaiseqinM0846dfgssgsd}
 \end{align}
 \end{lem}
\begin{pf}
 We begin with rewriting \eqref{01dsaf16asdfasf}$_2$ and \eqref{01dsaf16asdfasf}$_3$ into a nonhomogeneous form:
\begin{equation}\label{s0106pnnnn} \begin{cases}
 {\rho} u_t-\Delta(\mu u+\kappa\eta)+\nabla q =\mathcal{N}^1,\\[1mm]
\mathrm{div} {u} = -\mathrm{div}_{\tilde{\mathcal{A}}} {u},
\end{cases}  \end{equation}
where $\mathcal{N}^1:= \mu ( \mm{div}_{\tilde{\ml{A}}}\nabla_{\ml{A}}u + \mm{div} \nabla_{\tilde{\ml{A}}}u)
-\nabla_{\tilde{\ml{A}}}q $.

An application of $\partial^\alpha$ with $|\alpha |=2$ to \eqref{s0106pnnnn} yields
\begin{equation}   \label{201904291605xx}
  \begin{cases}
 {\rho}\partial^\alpha u_t+\nabla\partial^\alpha q-\Delta\partial^\alpha( \mu u+{\kappa} \eta)
=\partial^\alpha\mathcal{N}^1,   \\
 \mathrm{div} \partial^\alpha {u} =-\partial^\alpha\mathrm{div}_{\tilde{\mathcal{A}}} {u}.
  \end{cases}
\end{equation}
%% where $\alpha$ satisfies $|\alpha|=2$.
Then, multiplying \eqref{201904291605xx}$_1$ by $\partial^\alpha u$ in $L^2$, integrating by parts and using
\eqref{201904291605xx}$_2$, we get
\begin{align}
& \frac{1}{2}\frac{\mm{d}}{\mm{d}t}\left(    {\rho} \|\partial^\alpha u\|_0^2 + \kappa
\|\nabla  \partial^\alpha\eta\|_0^2 \right) +   \mu \|\nabla \partial^\alpha u \|_0^2 \nonumber \\
&=  \int  \partial^\alpha {\mathcal{N}}^1  \cdot   \partial^\alpha  u\mm{d}y  -\int \partial^\alpha  q   \partial^\alpha\mathrm{div}_{\tilde{\mathcal{A}}} {u} \mm{d}y   =:I_1+I_2 , \qquad (\;|\alpha|=2\;) \label{201808100154628xxx}
\end{align}
where the integral $I_1$ can be bounded as follows, using integration by parts, \eqref{20160614fdsa1957x} and \eqref{2017020614181721}.
\begin{align}
I_1& \lesssim  \| \nabla{\mathcal{N}}^1\|_0\|\nabla^3 u\|_0 \nonumber \\
&\lesssim( \|  \tilde{\mathcal{A}}  \|_{L^\infty}( (1+\|\mathcal{A}\|_{L^\infty})\| \nabla^3 u\|_0+\|\nabla^2 q\|_0   )+ \|\nabla  \tilde{\mathcal{A}} \|_{L^3}(
 (1+\|\tilde{\mathcal{A}}\|_{L^\infty})\|\nabla^2 u\|_{L^6} \nonumber \\
 &\quad
+\|\nabla  q\|_{L^6} )  +( \|\nabla^2 \mathcal{A}\|_{0}(1+\|\tilde{\mathcal{A}}\|_{L^\infty})+\|\nabla\tilde{\mathcal{A}}\|_{L^6}\|\nabla \mathcal{A}\|_{L^3})\|\nabla u\|_{L^\infty} )
\|\nabla^3 u\|_0\nonumber  \\
&\lesssim
\|\nabla \eta\|_2(\|\nabla u\|_2^2+
\| \nabla^3 u\|_0\|\nabla^2 q\|_0)  \nonumber \\ &\lesssim
  \|\nabla \eta\|_2\|\nabla^3  u\|_0^2  (1+\|\nabla u\|_1)+\kappa \| \nabla \eta\|_2^3\|\nabla^3 u\|_0. \label{201909122023}
\end{align}
%where we have used the integral by parts in the first inequality, and \eqref{20160614fdsa1957x} and \eqref{2017020614181721} in the last inequality.

Similarly, $I_2$ can be estimated as follows.
\begin{align}
I_2&\lesssim
\|\nabla^2q\|_0(\| \nabla^2\tilde{\mathcal{A}}\|_0 \|\nabla u\|_{L^\infty} +\|\nabla\tilde{\mathcal{A}}\|_{L^3}\|\nabla^2 u\|_{L^6}+\|\tilde{\mathcal{A}}\|_{L^\infty}\|\nabla^3 u\|_0) \nonumber \\
& \lesssim\|\nabla \eta\|_2 \|\nabla^3  u\|_0^2\|\nabla u\|_1 +\kappa\|\nabla\eta\|_2^3\|\nabla^3 u\|_0. \label{201909171446}
\end{align}

Plugging \eqref{201909122023} and \eqref{201909171446} into \eqref{201808100154628xxx}, and employing \eqref{aprpiosesnewxxxx},
\eqref{aprpiose1snewxxxxz} and Young's inequality, we obtain \eqref{ssebdaiseqinM0846dfgssgsd}.
\hfill$\Box$
\end{pf}
%%%%%%%%%%%%%%%%%%%%%%%%%%%%%%%%%%%%%%%%%%%%%%%%%%%%%%%%%
\begin{lem}\label{201612132242nxsfssdfsxx}
  Under the assumptions \eqref{aprpiosesnewxxxx}--\eqref{aprpiose1snewxxxxz} with sufficiently small $\delta$, it hold that
\begin{align}
\frac{\mm{d}}{\mm{d}t}\left(  {\rho} \sum_{|\alpha|=2} \int \partial^\alpha \eta \partial^\alpha u\mm{d}y +\frac{\mu}{2}\|\nabla^3 \eta\|_0^2\right)+
 \frac{\kappa}{2}  \|\nabla^3 \eta \|_0^2 \leqslant \rho  \| \nabla^2 u \|_{ 0}^2+ c\delta \|\nabla^3 u\|_0^2. \label{Lem:03xx}
 \end{align}
%% where $\alpha$ satisfies $ |\alpha|= 2$.
\end{lem}
\begin{pf}
We multiply \eqref{201904291605xx}$_1$  by $\partial^\alpha\eta$ ($|\alpha|=2$) in $L^2$ and use \eqref{01dsaf16asdfasf}$_1$ to get
\begin{align}
 \frac{\mm{d}}{\mm{d}t}  \left(  {\rho} \int \partial^\alpha \eta \partial^\alpha u\mm{d}y +\frac{\mu}{2}\|\nabla \partial^\alpha \eta \|_0^2\right)+
 {\kappa} \|\nabla \partial^\alpha \eta \|_0^2 =\rho\|\partial^\alpha u\|_0^2+ I_3+I_4, \label{2010154628}
\end{align}
where
$$I_3:=\int\partial^\alpha {\mathcal{N}}^1\cdot\partial^\alpha\eta\mm{d}y \;\;
\mbox{ and }\;\; I_4:=\int \partial^\alpha q\partial^\alpha\mathrm{div}\eta \mm{d}y. $$

Similarly to \eqref{201909122023} and \eqref{201909171446}, $I_3$ and $I_4$ can be controlled as follows.
\begin{align}   \label{201911262125}
|I_3|\lesssim &  \| \nabla{\mathcal{N}}^1\|_0\| \nabla^3\eta\|_{0}   \lesssim
\|\nabla \eta\|_2^2\|\nabla^3  u\|_0  (1+\|\nabla u\|_1) + \kappa \| \nabla \eta\|_2^4     \nonumber
\end{align}
and
\begin{align}
|I_4|\lesssim \|\mm{div}\eta\|_2\|\nabla^2 q\|_0  \lesssim  \|\nabla \eta\|_2^2
    (\|\nabla u\|_1^{2}+ {\kappa}  \| \nabla \eta\|_2^2),
\end{align}
where we have used \eqref{2016957x}  in the second inequality in \eqref{201911262125}.

Putting the above estimates into \eqref{2010154628}, and using \eqref{aprpiosesnewxxxx}--\eqref{aprpiose1snewxxxxz}
and \eqref{20160614fdsa1957x}, we get \eqref{Lem:03xx}. \hfill$\Box$
 \end{pf}

Now, we are in a position to show \eqref{201905041053}. If one makes use of \eqref{20160614fdsa1957x},
\eqref{ssebdaiseqinM0846dfgssgsd} and \eqref{Lem:03xx}, then one has
\begin{equation}
\label{201911262012}
\frac{\mm{d}}{\mm{d}t}  {\mathcal{E}}_1+
c\|\nabla( u,\sqrt{\kappa}\eta)\|_2^2 \leqslant\rho\|\nabla^2 u\|_{0}^2,
\end{equation}
 where
$$ {\mathcal{E}}_1: = {c}_4({\rho} \| \nabla^2 u\|_0^2+ \kappa \| \nabla^3 \eta \|_0^2) +
{\rho} \sum_{  |\alpha| = 2}\int \partial^\alpha \eta \partial^\alpha u\mm{d}y +\frac{\mu}{2}\|\nabla^3 \eta\|_0^2$$
satisfies $\|\nabla u\|_1^2 +\|\sqrt{\kappa}\nabla\eta\|_2^2\lesssim {\mathcal{E}}_1$.

Noting that $(\eta,u)$ satisfies \eqref{ssebdaiseqinM0846dfgssgsdxx},
using \eqref{20160614fdsa1957}, \eqref{20160614fdsa1957x} and the interpolation inequality \eqref{201905041130},
we further get from \eqref{ssebdaiseqinM0846dfgssgsdxx} and \eqref{201911262012} that there are constants $\delta_1\in (0,1]$,
${c}_5$ and ${c}$, such that
\begin{equation}    \label{201811262026}
\frac{\mm{d}}{\mm{d}t} {\mathcal{E}}_2 + c(\| \nabla ( u,\sqrt{\kappa}\eta) \|_2^2 )  \leqslant 0
\end{equation}
for any $K$ and $\kappa$ satisfying \eqref{aprpiosesnewxxxx} and \eqref{aprpiose1snewxxxxz} with $\delta\leqslant\delta_1$, where  ${\mathcal{E}}_2:={\mathcal{E}}_1 + {c}_5(\rho\|u\|_0^2 +\kappa\|\nabla\eta\|_0^2)$ satisfies
\begin{equation}\nonumber
\| (u,\sqrt{\kappa}\nabla \eta)\|_2^2  \lesssim  {\mathcal{E}}_2 \lesssim\|( u,\sqrt{\kappa}\nabla \eta)\|_2^2  .
\end{equation}

An integration of \eqref{201811262026} over $(0,t)$ yields that there is a constant ${c}_1\geqslant 1$, such that
\begin{equation}
\label{201811262026xx}
\|  ( u,\sqrt{\kappa}\nabla \eta)\|_2^2 +c  \int_0^t\| \nabla ( u,\sqrt{\kappa}\eta)\|_2^2 \mm{d}\tau
\leqslant {c}_1{I_0^{\mm{h}}( u^0, \eta^0)}\quad\mbox{for any }t\in [0,T),
\end{equation}
which gives the desired stability estimate \eqref{201905041053}. Finally, if we take
\begin{equation}
\label{201911262060}
K:= 2\sqrt{ {c}_1 {I_0^{\mm{h}}(  u^0, \eta^0)}} ,
\end{equation}
 we obtain \eqref{aprpiose1} under the assumptions \eqref{aprpiosesnewxxxx} and \eqref{aprpiose1snewxxxxz} with $\delta\leqslant\delta_1$.

Next we introduce the local well-posedness for  the initial value problem \eqref{01dsaf16asdfasf}--\eqref{01dsaf16asdfasfsaf}.
%%%%%%%%%%%%%%%%%%%%%%%%%%%%%%%%%%%%%%%
\begin{pro} \label{pro:0401nxdxx}
Let $B_1\geqslant 0$, $(\eta^0,u^0)\in H^3 \times H^2$ satisfy $\|u^0\|_2\leqslant B_1$, let $\zeta^0:= \eta^0+y$ and $\mathcal{A}^0$ be
defined by $\zeta^0$. Then there is a constant $\delta_2\in (0,1]$, such that for $(\eta^0,u^0)$
satisfying $\mm{div}_{\mathcal{A}^0}u^0=0$ and
\begin{align}
  \|\nabla \eta^0\|_2 \leqslant \delta_2,    \label{201912261426}
\end{align}
the initial value problem \eqref{01dsaf16asdfasf}--\eqref{01dsaf16asdfasfsaf} has a unique strong solution
$( \eta, u,q)\in C^0([0,T ),H^3_{*,1})\times \mathfrak{U}_T\times L^\infty(I_T, \underline{H}^2 )$ for some $T>0$ which may depend
on $B_1$, $\rho$ and $\mu$.
\end{pro}
\begin{pf}
The proof of Proposition \ref{pro:0401nxdxx} will be provided in Section \ref{201912062021}.
  \hfill $\Box$
\end{pf}

With the \emph{priori} estimate \eqref{201811262026xx} (under the assumptions \eqref{aprpiosesnewxxxx} and \eqref{aprpiose1snewxxxxz} with $\delta\leqslant\delta_1$) and Proposition \ref{pro:0401nxdxx} in hand, we can easily establish Theorem \ref{201904301948} as follows.
Next, we briefly give the proof.

{\sc Proof of Theorem \ref{201904301948}}. \  Let $(\eta^0,u^0)\in H^3_{*,1} \times H^2$, and $\kappa$ be properly large so that
\begin{align}
 \max\{K , K^4\}/\kappa\leqslant\min\{\delta_1, \delta_2^2/2, 2{  {c}_1\delta_2^2}\},
  \label{201912261425}
  \end{align}
  where $K$ is defined by \eqref{201911262060}.

Noting that
\begin{align}  \label{202001191736}
K^2\leqslant \frac{2K}{3}+\frac{K^4}{3}\leqslant 2\max\{K , K^4\} ,
\end{align}
we see that $\eta^0$ satisfies \eqref{201912261426} due to \eqref{201912261425}. Thus, by virtue of Proposition \ref{pro:0401nxdxx},
there exists a unique local solution
$(\eta, u,q)\in C^0([0,T^{\max}),H^3_{*,1})\times \mathfrak{U}_{T^{\max}}\times L^\infty(I_{T^{\max}},\underline{H}^2)$
to the initial value problem \eqref{01dsaf16asdfasf}--\eqref{01dsaf16asdfasfsaf}, where $T^{\max}$ is the maximal time of existence.
Let
\begin{equation}   \nonumber
T^{*}:=\sup\left\{ T \in I_{T^{\max}}~\left|~ \|(u(t),\sqrt{\kappa}\nabla\eta(t))\|_2\leqslant K\mbox{ for any }t\leqslant T\right.\right\}.
\end{equation}
Recalling the definition of $K$, we find that the definition of $T^*$ makes sense and $T^*>0$.

Now, we prove $T^*=\infty$ by contradiction. Assume that $T^*<\infty$. Noting that $T^{\max}$ denotes the maximal existence time and $K/\sqrt{\kappa}\leqslant \delta_2$ by virtue of \eqref{201912261425} and \eqref{202001191736}, we apply Proposition \ref{pro:0401nxdxx}
and recall the definition of $T^{\max}$ to easily see that
\begin{equation}
\label{201911262202}
\limsup_{t\to T^*}\|(u(t),\sqrt{\kappa}\nabla \eta(t))\|_2^2=K.
\end{equation}

Since $\max\{K,K^4\}/\kappa\leqslant\delta_1$ and $\sup_{0\leqslant t< T^*}\|( u(t),\sqrt{\kappa}\nabla \eta(t))\|_2\leqslant K$,
 we can show that the solution $(\eta,u )$ enjoys the stability estimate \eqref{201811262026xx} with $T^*$ in place of $T$ by the regularity
 of $(\eta,u,q)$. More precisely, we have
$$\|(u(t),\sqrt{\kappa} \nabla \eta(t)\|_2\leqslant \sqrt{ {c}_1} \|(u^0,\sqrt{\kappa}\nabla \eta^0)\|_2\leqslant K/ {2}
\;\;\mbox{ for any }t< T^*,$$
which contradicts with \eqref{201911262202}. Hence, $ T^*=\infty$.

This completes the proof of the existence and the estimate \eqref{201905041053} in Theorem \ref{201904301948} with
$ {c}_2:=\min\{\delta_1, \delta_2^2/2, 2  {c}_1\delta_2^2\}$. The uniqueness of the global solution is obvious
due to the uniqueness of the local solution. In addition, \eqref{202001121649}--\eqref{202001121644} are obvious to obtain.

\section{Proof of Theorem \ref{201912041028}} \label{201912042253}

This section is devoted to the proof of Theorem \ref{201912041028}. We start with the proof of the first conclusion.
\subsection{Exponential stability}\label{202001192326}
Let $\bar{u}=u-({u}^0)_{\mathbb{T}^3}$ and $\bar{u}^0:=u^0-( {u}^0)_{\mathbb{T}^3}$, and we rewrite \eqref{01dsaf16asdfasf} as follows.
\begin{equation}
                              \begin{cases}
\eta_t =\bar{u}+({u}^0)_{\mathbb{T}^3} , \\[1mm]
 {\rho} \bar{u}_t-\mu \Delta_{\ml{A}}\bar{u}+\nabla_{\ml{A}}q=
{\kappa}\Delta \eta, \\[1mm]
\div_\ml{A}\bar{u}=0
\end{cases} \label{202003031102}
\end{equation}
with initial value  condition
\begin{equation}\nonumber  (\eta, \bar{u})|_{t=0}=(\eta^0,  \bar{u}^0)\mbox{ in }\mathbb{T}^3.
\end{equation}
Then, following the arguments used in the derivation of \eqref{201811262026}, we have
\begin{equation}
\label{2018112620261x}
\frac{\mm{d}}{\mm{d}t}  \bar{\mathcal{E}}_2+
c\| \nabla (\bar{u},\sqrt{\kappa}\eta) \|_2^2   \leqslant 0,
\end{equation}
where $\bar{\mathcal{E}}_2$ is defined as ${\mathcal{E}}_2$ with $u$ replaced by $\bar{u}$, and satisfies
\begin{equation}
 \label{209120122308ss} \|(\bar{u},\sqrt{\kappa}\nabla \eta )\|_2^2 \lesssim \bar{\mathcal{E}}_2 \lesssim  \|(\bar{u},\sqrt{\kappa}\nabla \eta )\|_2^2.
\end{equation}
Noting that $(\bar{u})_{\mathbb{T}^3}=0$, we further deduce from \eqref{2018112620261x} and  \eqref{209120122308ss} that
\begin{equation}
\label{2018df11262026xx}
\frac{\mm{d}}{\mm{d}t}  \bar{\mathcal{E}}_2+  {c}_3 ( \bar{\mathcal{E}}_2+ \|\bar{u}\|_3^2+\kappa\| \nabla  \eta\|_2^2)   \leqslant 0,
\end{equation}
which implies \eqref{201905041x053} immediately.

\subsection{Large-time behavior of $\eta$}
Let $\bar{\eta}:=\eta-\eta(y_{\mm{h}},0 ,t)$, then by virtue of \eqref{201905041x053},
 \begin{align}
  \|\partial_3\bar{\eta}\|_0^2 +\int_0^t \|\partial_3\bar{\eta}\|_0^2  e^{ {c}_3(\tau-t)}\mm{d}\tau
  \lesssim e^{- {c}_3t}I_0^{\mm{h}}(\bar{u}^0,\eta^0)/{\kappa}. \nonumber
\end{align}
Recalling $ \bar{\eta}|_{y_3=0}=\bar{\eta}|_{y_3=2\pi}=0$, we use Wirtinger's inequality (see (4.27) in \cite{JFJSOMITNx}
or the first formula after (3.56) in \cite{JFJSOUI}) to infer that
 \begin{align}
   \| \bar{\eta}\|_0^2 +\int_0^t \| \bar{\eta} \|_0^2  e^{ {c}_3(\tau-t)}\mm{d}\tau
   \lesssim e^{- {c}_3t}I_0^{\mm{h}} (\bar{u}^0,\eta^0) /{\kappa} .\label{201905041x053x}
\end{align}

In view of the trace estimate (see \cite[Lemma 9.7]{JFJSZWCIARTPA})
\begin{align}
 \|f|_{y_3=a}\|_{H^{1/2}((0,2\pi)^2)}\lesssim_0\|f\|_{1}\;\;\mbox{ for any }f\in H^{1} \mbox{ and for any } a\in (0,2\pi), \nonumber
\end{align}
we have
$$ \begin{aligned}
\|\nabla \bar{\eta}\|_1^2\leqslant & 2(  \|\nabla {\eta}\|_1^2 +\|\nabla \eta(y_{\mm{h}},0 ,t) \|_1^2)\\
=&2\left( \|\nabla {\eta}\|_1^2  +\sum_{1\leqslant \alpha_1+\alpha_2\leqslant 2}\|\partial_1^{\alpha_1}\partial_2^{\alpha_2}
\eta(y_{\mm{h}},0 ,t)\|_0^2\right)\lesssim \|\nabla \eta \|_2^2.
\end{aligned}$$
Thus, by \eqref{201905041x053},
$$   \begin{aligned}
 & \|\nabla \bar{\eta}\|_1^2+\int_0^t \|\nabla \bar{\eta} \|_1^2  e^{ {{c}}_3(\tau-t)}\mm{d}\tau   \lesssim e^{- {c}_3t}I_0^{\mm{h}}(\bar{u}^0,\eta^0) /\kappa ,
 \end{aligned}
 $$
 which, together with \eqref{201905041x053x}, yields \eqref{201905041053x}.

Now we turn to derive \eqref{20191211197}. By \eqref{202003031102}$_1$, we have
$$\eta-({u}^0)_{\mathbb{T}^3}t=f(y,t):=\eta^0+\int_0^t\bar{u}\mm{d}\tau. $$
From the estimate of $\nabla \eta$ in \eqref{201905041x053} we easily get
\begin{align}
\nabla f(y,t)\to 0 \mbox{ strongly in }H^2\mbox{ as }t\to \infty.
\label{202003031128}
\end{align}

We use the estimate of $\bar{u}$ in \eqref{201905041x053} to deduce that
$$\|f(y,t)\|_2^2\lesssim \|\eta^0\|_2^2+I_0^{\mm{h}}(\bar{u}^0 , \eta^0). $$
Therefore, there exists a sequence $\{t_n\}_{n=1}^\infty$, such that
$$f(y,t_n)\to \varpi\mbox{ strongly in }H^1\mbox{ as }t_n\to \infty. $$
In view of \eqref{202003031128}, we immediately see that $\varpi$ is a constant vector.

Recalling the definition of $f$ and the fact $(\bar{u})_{\mathbb{T}^3}=0$, one has
$$\int \eta^0\mm{d}y=\int f(y,t_n)\mm{d}y \to \int \varpi\mm{d}y, $$
which implies $\varpi=(\eta^0)_{\mathbb{T}^3}/8\pi^3$. Since the vector $\varpi$ is uniquely determined by $\eta^0$,
we can easily prove by contradiction that
$$f(y,t)\to \varpi\mbox{ strongly in }H^1\mbox{ as }t\to \infty. $$
  Thus, for any $t\geqslant 0$,
\begin{align}
\label{202003031153}
(f(y,t)-\varpi)_{\mathbb{T}^3}=0.
\end{align}
Thanks to \eqref{201905041x053},  \eqref{202003031153} and Poinc\'are's inequality, we obtain \eqref{20191211197}.

\subsection{Stability of $(\eta,u)$ around $(\eta^{\mm{l}},u^{\mm{l}})$}   \label{202001081357}

 Let $(\eta^0,u^0)$ be given in Theorem \ref{201904301948}.
 Then by the regularity theory of the Stokes problem, there exists a unique solution $(\eta^{\mm{r}},u^{\mm{r}},Q_1,Q_2)$, satisfying
\begin{align}
\label{202001072004xx}
   \begin{cases}
   -\Delta \eta^{\mm{r}} +\nabla Q_1=0, \\
   \mm{div} \eta^{\mm{r}}=-\mm{div}\eta^0,\\
     (\eta^{\mm{r}})_{\mathbb{T}^3}  = 0
   \end{cases}
\end{align}
and
\begin{align}
\label{202001072004}
   \begin{cases}
   -\Delta u^{\mm{r}} +\nabla Q_2=0, \\
   \mm{div} u^{\mm{r}}=\mm{div}_{\tilde{\mathcal{A}}^0}u^0,\\
     ( {u}^{\mm{r}})_{\mathbb{T}^3}  = 0,
   \end{cases}
\end{align}
where $\tilde{\mathcal{A}}^0:={\mathcal{A}}^0-I$. Moreover,
\begin{align}
&\|\eta^{\mm{r}}\|_3\lesssim \|\mm{div}\eta^0\|_2\lesssim  \|\nabla \eta^0\|_2^2 ,
\label{202001121445x} \\
&\|u^{\mm{r}}\|_2\lesssim \|\mm{div}_{\tilde{\mathcal{A}}^0}u^0\|_1\lesssim \|\tilde{\mathcal{A}}^0\|_2\|u^0\|_2\lesssim \|\nabla \eta^0\|_2\|u^0\|_2, \label{202001121445}
\end{align}
where we have used the fact $\|\nabla \eta^0\|_2\lesssim_0 1$  by \eqref{202001121644} in the estimate \eqref{202001121445x}.

Let $\tilde{\eta}^0=\eta^0+\eta^{\mm{r}}$ and $\tilde{u}^0=u^0+u^{\mm{r}}$. Then it is easy to see that
$(\tilde{\eta}^0,\tilde{u}^0)\in H^3_\sigma\times H^2_\sigma$. Thus there exists a unique global solution
$(\eta^{\mm{l}},u^{\mm{l}} )\in  C^0(\mathbb{R}_0^+, H^{3})\times \mathfrak{U}_\infty $  to \eqref{01dsaf16asdfasfxx}
with initial condition $(\eta^{\mm{l}}, {u}^{\mm{l}})|_{t=0}=(\tilde{\eta}^0,\tilde{u}^0)$ and with
\begin{align}
 \nabla q^{\mm{l}}=0\label{1201905asfda041x053xx}.
\end{align}
 For the proof of the above claim, we refer to the proof of Proposition \ref{qwepro:0sadfa401nxdxx} in Section \ref{2020010813027}.
Moreover,
\begin{align}
\label{202001121635}
\mm{div}\eta^{\mm{l}}=\mm{div}\tilde{\eta}^0=0.
 \end{align}

Let
$\bar{u}^{\mm{l}}:=u^{\mm{l}}-(\tilde{u}^0)_{\mathbb{T}^3}$ and  $\bar{u}^{\mm{l}}_0:=u^0+u^{\mm{r}}-( {u}^0)_{\mathbb{T}^3}$. Then,
similarly to \eqref{201905041x053}, and utilizing \eqref{202001121445x} and \eqref{202001121445}, one can also get
\footnote{ It should be remarked that three constants ${c}_3$ in \eqref{2018df11262026xx}, \eqref{201905asfda041x053}
and \eqref{201811sdf262026} can be taken to be equal.}
\begin{align}
&  \|(  \bar{u}^{\mm{l}} ,\sqrt{\kappa}\nabla \eta^{\mm{l}})\|_2^2 +\int_0^t( \| \bar{u}^{\mm{l}} \|_3^2+{\kappa}\|\nabla \eta^{\mm{l}}\|_2^2 )e^{ {c}_3(\tau-t)}\mm{d}\tau\nonumber \\
& \qquad \lesssim e^{- {c}_3 t} I_0^{\mm{h}}( \bar{u}^{\mm{l}}_0, \tilde{\eta}^0)\lesssim e^{- {c}_3 t}
( \| \nabla \eta^0 \|_2^2 I_0^{\mm{h}}( {u}^0, {\eta}^0) +I_0^{\mm{h}} (\bar{u}^0, {\eta}^0)).  \label{201905asfda041x053}
\end{align}
%% where we have used \eqref{202001121445x} and \eqref{202001121445} in the last inequality in \eqref{201905asfda041x053}.

Now, if we subtract \eqref{01dsaf16asdfasfxx} and \eqref{202001070914} from the initial value problem \eqref{01dsaf16asdfasf} and \eqref{01dsaf16asdfasfsaf}, then we get
 \begin{equation}\label{01dsaf16safafasdfasfxx}
                              \begin{cases}
\eta_t^{\mm{d}}=u^{\mm{d}}, \\[1mm]
 {\rho} u_t^{\mm{d}}+\nabla q^{\mm{d}}-\Delta ( \mu u^{\mm{d}}+
{\kappa} \eta^{\mm{d}})= \mathcal{N}^1, \\[1mm]
\div u^{\mm{d}}=-\mathrm{div}_{\tilde{\mathcal{A}}} {u},\\
(\eta^{\mm{d}},u^{\mm{d}})|_{t=0}=-(\eta^{\mm{r}},u^{\mm{r}}).
\end{cases}
\end{equation}
Thus we can follow the arguments used in Lemmas \ref{201612132242nx}--\ref{201612132242nxsfssdfsxx} with necessary slight modifications
to deduce from \eqref{01dsaf16safafasdfasfxx} that
\begin{align}
 & \frac{\mm{d}}{\mm{d}t}  \left(  {\rho} \| \nabla^i u^{\mm{d}}\|_0^2+ \kappa \| \nabla^{i+1} \eta^{\mm{d}} \|_0^2\right)+
 \mu \|\nabla^{i+1}  u^{\mm{d}} \|_0^2\nonumber \\
 & \qquad \lesssim \|\nabla \eta\|_2\|\nabla^i u\|_1\|\nabla q^{\mm{d}} \|_1
+\|\nabla^{i} u^{\mm{d}} \|_1\|\mathcal{N}^1\|_1 , \label{201702061safa4181721asfdx} \\
&  \frac{\mm{d}}{\mm{d}t}\Big(  {\rho} \sum_{|\alpha|=2} \int \partial^\alpha \eta^{\mm{d}} \partial^\alpha u^{\mm{d}}\mm{d}y +\frac{\mu}{2}\|\nabla^3 \eta^{\mm{d}}\|_0^2\Big)+  {\kappa}  \|\nabla^3\eta^{\mm{d}} \|_0^2  \nonumber \\
 & \qquad \leqslant\rho\|\nabla^2u^{\mm{d}}\|_0^2+ \|\mm{div} \eta^{\mm{d}}\|_2\|\nabla q^{\mm{d}} \|_1
 + \|\nabla \eta^{\mm{d}} \|_2\|\mathcal{N}^1\|_1 , \label{Lem:03sfdaxx}
 \end{align}
where $i=0$ and $2$.

In view of  \eqref{202001121644}, one has
\begin{align}
& \|\mathcal{N}^1\|_1\lesssim \|\nabla \eta\|_2(\|\nabla u \|_2+\| \nabla q\|_1 ).\label{201912270953}
\end{align}
In addition, we employ \eqref{201909261909} with $J=1$, \eqref{202001072004xx}$_3$ and \eqref{202001072004}$_3$ to deduce from \eqref{01dsaf16safafasdfasfxx}$_1$ and \eqref{01dsaf16safafasdfasfxx}$_2$ that
\begin{align}
(u^{\mm{d}})_{\mathbb{T}^3}=0  \mbox{ and }(\eta^{\mm{d}})_{\mathbb{T}^3}=0\mbox{ for any } t\geqslant 0. \label{201912271373}
\end{align}
Consequently, with the help of \eqref{1201905asfda041x053xx}, \eqref{202001121635}, \eqref{201912270953} and \eqref{201912271373},
we argue, similarly to \eqref{2018df11262026xx}, to get from \eqref{201702061safa4181721asfdx}--\eqref{Lem:03sfdaxx} that
\begin{align}
&\frac{\mm{d}}{\mm{d}t}  {\mathcal{E}}_2^{\mm{d}}+  {c}_3({\mathcal{E}}_2^{\mm{d}}
+\|u^{\mm{d}} \|_3^2 +\kappa \| \nabla \eta^{\mm{d}}\|_2^2)
 \lesssim \|\nabla q \|_1 (\|\nabla  \eta \|_2 ( \|u\|_1+\|\nabla   u \|_2)+\|\mm{div}\eta \|_2)  \nonumber \\
 & \qquad + \|\nabla \eta\|_2 \|\nabla (\eta^{\mm{d}},u^{\mm{d}}) \|_2 (\|\nabla u \|_2+\| \nabla q\|_1 ) ,
\label{201811sdf262026}
\end{align}
where ${\mathcal{E}}_2^{\mm{d}}$, defined as ${\mathcal{E}}_2$ with $(\eta ,u )$ replaced by $(\eta^{\mm{d}},u^{\mm{d}})$,
satisfies the estimate
  $$ \|  u^{\mm{d}}\|_2^2 +\kappa\| \eta^{\mm{d}}\|_3^2 \lesssim  {\mathcal{E}}_2^{\mm{d}}
  \lesssim \|  u^{\mm{d}}\|_2^2 +\kappa\| \eta^{\mm{d}}\|_3^2. $$

Multiplying \eqref{201811sdf262026} by $e^{ {c}_3 t}$, integrating the resulting inequality over $(0,t)$, we make use of \eqref{201905041053}--\eqref{202001071328}, and \eqref{201905041x053}, \eqref{201905asfda041x053} and Young's inequality
to obtain
 \begin{align}
& \| u^{\mm{d}}\|_2^2+ {\kappa}\| \eta^{\mm{d}}\|_3^2+  \int_0^t(\|u^{\mm{d}} \|_3^2
 +\kappa \| \nabla \eta^{\mm{d}}\|_2^2) e^{ {c}_3(\tau-t)} \mm{d}\tau\nonumber\\
 & \quad \lesssim \sqrt{\kappa^{-1}}\max\{1,\sqrt{\kappa}^{-1}\}e^{- {c}_3 t}\nonumber \\
&\qquad\; \times \left(\sqrt{I_0^{\mm{h}}(u^0,\eta^0)} +I_0^{\mm{h}}(u^0,\eta^0) \right)
 {( \| \nabla \eta^0 \|_2^2 I_0^{\mm{h}}(  {u}^0, {\eta}^0) +I_0^{\mm{h}} ( \bar{u}^0, {\eta}^0))}.\nonumber
\end{align}
Recalling $(\eta^{\mm{d}})_{\mathbb{T}^3}=0 $, we further conclude \eqref{201905041053xxx}.
This completes the proof of Theorem \ref{201912041028}.

\section{Proof of Theorems \ref{x2019saf04301948} and \ref{20191asfas2041028}} \label{20201011926}

In this section we give the proof of Theorems \ref{x2019saf04301948} and \ref{20191asfas2041028}.
Since the proof is similar to that of Theorems \ref{201904301948} and \ref{201912041028}, next we only sketch the proof of
Theorems \ref{x2019saf04301948} and \ref{20191asfas2041028}.

The key step to get Theorem \ref{x2019saf04301948} is also to establish the \emph{a priori} stability estimate. To this end,
let $(\eta,u)$ be a solution of the initial value problem \eqref{n0101nnn} and \eqref{2201912032006}, and satisfy
$\eta^0 \in H^{3}_{*}$,
\begin{align}
 \label{aprpiosesnewxxxxx}
 \sup_{0\leqslant t < T}\| (u(t),\sqrt{\kappa}\nabla \eta(t))\|_2 \leqslant K_{\mm{P}}\;\;\mbox{ for any given }T,
\end{align}
and
\begin{equation}\label{aprpiosesnewxxxxzx}
\max\{K_\mm{P}^{2/3} , K^2_\mm{P}\}/\kappa\leqslant \delta^2\in (0,1],\end{equation}
 where $K_{\mm{P}}$ will be given by \eqref{201saf911262060x} and $\delta$ is a sufficiently small constant.
 We remark here that the smallness of $\delta$ only depends on $P(\cdot)$, $\bar{\rho}$, $\mu$ and $\lambda$.

  In addition, under the assumptions \eqref{aprpiosesnewxxxxx} and \eqref{aprpiosesnewxxxxzx}, one has
\begin{equation}
 \sup_{0\leqslant t < T} \|\nabla \eta(t)\|_2   \leqslant  \delta .
\label{202001191952}
\end{equation}
Then, the following  preliminary estimates can be established.
\begin{lem}
\label{201805141072dsafa} Let $\eta$ satisfy \eqref{202001191952} with sufficiently small $\delta $, then
\begin{align}
& 1/2 \leqslant J\leqslant 3/2,  \label{201912302317} \\
& \|J-1\|_2\lesssim_0  \|\nabla  \eta\|_2,\label{20saf1912302317} \\
&\|J^{-1}-1+\mm{div}\eta\|_2\lesssim_0  \|\nabla \eta\|_2^2,  \label{2019123023171}  \\
&\|\tilde{\mathcal{A}}\|_2\lesssim_0   \|\nabla \eta\|_2  ,\label{aimasdadsdfse} \end{align}
where $J:=|\det(\nabla \eta+I)|$.
\end{lem}
%%%%%%%%%%%%%%%%%%%%%%%%%%%%%%%%%%%%
\begin{pf}
Recalling the second identity in \eqref{202001081820} and the derivation of \eqref{2017020614181721xxdfasdfasadfaff}, we see that
\begin{equation}
\label{Jdetmentionsa}
J=1+  \mm{div} \eta+r_\eta
\end{equation}
and
\begin{equation}
\label{Jdetmentsadfafionsa}
\|r_\eta\|_2 \lesssim_0    \| \nabla \eta\|_2^2.
\end{equation}
Then, from \eqref{202001191952}, \eqref{Jdetmentionsa} and  \eqref{Jdetmentsadfafionsa}  we get \eqref{20saf1912302317}, while
  the estimate \eqref{201912302317} follows from \eqref{201905041149} and \eqref{20saf1912302317}.

In view of \eqref{Jdetmentionsa}, we find that
\begin{equation}J^{-1}-1+
 \mm{div}\eta=(1-J^{-1}- \mm{div}\eta )\left(\mm{div}\eta+r_\eta\right)
 + \mm{div}\eta\left(\mm{div}\eta+r_\eta\right) -r_\eta .\nonumber \end{equation}
  Thus we can easily  derive \eqref{2019123023171}   from  the above relation.

By \eqref{2019123023171} we see that
\begin{align}
\|\nabla J^{-1} \|_1\lesssim_0  \|\nabla \eta\|_2. \nonumber
\end{align}
Recalling the derivation of \eqref{202001192046}, we have
\begin{align}
 \|(A^*_{ij})_{3\times 3}-I\|_2 \lesssim_0  \| \nabla \eta\|_2^2. \nonumber
   \end{align}
Thanks to \eqref{201912302317},  \eqref{20saf1912302317} and the above two estimates, we easily obtain \eqref{aimasdadsdfse}
from the following relation
$$\tilde{\mathcal{A}}=J^{-1}({(A^*_{ij})_{3\times 3}^{\mm{T}}-I}+(1-J)I).$$
This completes the proof of Lemma \ref{201805141072dsafa}. \hfill$\Box$
\end{pf}

Next, we continue to derive more estimates for $(\eta,u)$.  %% as in Lemmas \ref{201612132242nn} and \ref{201612132242nxsfssdfsxx}.
\begin{lem}\label{qwe201612132242nn}
Under the conditions \eqref{aprpiosesnewxxxxx}--\eqref{aprpiosesnewxxxxzx} with sufficiently small $\delta$, we have
\begin{align} &
\frac{1}{2}\frac{\mm{d}}{\mm{d}t}  \Big(  \bar{\rho} \| u\|_0^2+ {\bar{\rho}}\int
\int_{\bar{\rho}/4}^{\bar{\rho}J^{-1}}\frac{P(z)}{z^2}\mm{d}z\mm{d}y + \kappa \| \nabla   \eta \|_0^2\Big)
+  \mu \|\nabla_{\mathcal{A}} u \|_0^2  + \lambda  \|\mm{div}_{\mathcal{A}}u\|_0^2 \leqslant 0,  \label{201912302325} \\
&  \frac{\mm{d}}{\mm{d}t}  \left(  \bar{\rho} \|\nabla^2  u\|_0^2 +   P'(\bar{\rho})\bar{\rho}\|\nabla^2 \mm{div}\eta\|^2_0
 + \kappa \| \nabla^3  \eta \|_0^2 \right)
 +  \mu \|\nabla^3  u \|_0^2  + \lambda  \|\nabla^2 \mm{div}  u\|_0^2 \lesssim_{\mm{P}}
\delta \kappa \|\nabla \eta\|_2^2.  \label{qwessebdaiseqinM0846dfgssgsd}
 \end{align}
 \end{lem}
\begin{pf}
(1) Multiplying \eqref{n0101nnn}$_2$ by $u$ in $L^2$ and using \eqref{201909261909}, we obtain
\begin{equation}
\frac{1}{2}\frac{\mm{d}}{\mm{d}t} \Big(  \bar{\rho} \| u\|_0^2+ 2\int J\nabla_{\mathcal{A}}P(\bar{\rho}J^{-1})\cdot u\mm{d}y
+\kappa\|\nabla \eta\|_0^2\Big)
+  \mu \|\sqrt{J}\nabla_{\mathcal{A}} u \|_0^2  + \lambda  \|\sqrt{J}\mm{div}_{\mathcal{A}}u\|_0^2 \leqslant 0.\label{201912302324}
\end{equation}

By virtue of \eqref{Jtdimau},
$$  \int J\nabla_{\mathcal{A}}P(\bar{\rho}J^{-1})\cdot u\mm{d}y= {\bar{\rho}}\frac{\mm{d}}{\mm{d}t}\int
\int_{\bar{\rho}/4}^{\bar{\rho}J^{-1}}\frac{P(z)}{z^2}\mm{d}z\mm{d}y. $$
Substituting the above identity into \eqref{201912302324} and utilizing \eqref{201912302317}, we get \eqref{201912302325}.

(2) Now we turn to the proof of \eqref{qwessebdaiseqinM0846dfgssgsd}. We begin with rewriting \eqref{n0101nnn}$_2$ into
a nonhomogeneous form:
\begin{align}
&\bar{\rho} u_t-\nabla( P'(\bar{\rho})\bar{\rho}\,\mm{div}\eta+\lambda\mm{div}u )-\Delta(\mu u +{\kappa}\eta)=\mathcal{N}^2,
\label{qwes0106pnnnn}\end{align}
where
$$
\begin{aligned}
 \mathcal{N}^2:= &\mu  ( \mm{div}_{\tilde{\ml{A}}}\nabla_{\ml{A}}u
+ \mm{div} \nabla_{\tilde{\ml{A}}}u)+\lambda  ( \nabla_{\tilde{\ml{A}}}\mm{div}_{\ml{A}}u
+ \nabla\mm{div}_{\tilde{\ml{A}}}u)\\
& + (J-1)(\mu \Delta_{\mathcal{A}}u+\lambda \nabla_{\mathcal{A}}\mm{div}_{\mathcal{A}}u)-(J-1) \nabla_{\mathcal{A}}
P(\bar{\rho}J^{-1})-\nabla_{\tilde{\mathcal{A}}}  P(\bar{\rho}J^{-1})  \\
& -\nabla \left( {P}'(\bar{\rho})\bar{\rho}(J^{-1}-1+\mm{div}\eta)  +\int_{0}^{\bar{\rho}(J^{-1}-1)}(\bar{\rho}
(J^{-1}-1)-z)\frac{\mm{d}^2}{\mm{d}z^2} P (\bar{\rho}+z)\mm{d}z\right) .
\end{aligned}
$$

Taking $\alpha$ with $|\alpha|=2$ and applying $\partial^{\alpha}$ to \eqref{qwes0106pnnnn}, then
multiplying \eqref{qwes0106pnnnn} by $\partial^{\alpha} u$ in $L^2$, we find that
\begin{align}
& \frac{1}{2}\frac{\mm{d}}{\mm{d}t}\left(    \bar{\rho} \|\partial^\alpha u\|_0^2
+  P'(\bar{\rho})\bar{\rho}\| \partial^\alpha\mm{div}\eta\|^2_0 + \kappa
\|\nabla  \partial^\alpha\eta\|_0^2 \right)
\nonumber \\
&+
  \mu \|\nabla \partial^\alpha u \|_0^2 + \lambda  \|\mm{div\partial^\alpha}u\|_0^2=   \int   \partial^\alpha {\mathcal{N}}^1  \cdot   \partial^\alpha  u\mm{d}y  =:I_5 . \label{qwe201808100154628xxx}
\end{align}
Thanks to Lemma \ref{201805141072dsafa}, we can follow the proof procedure used for \eqref{201909122023} to deduce that
\begin{align}
I_5\leqslant \| \nabla{\mathcal{N}}^2\|_0\|\nabla^3 u\|_0 \lesssim_{\mm{P}}
\|\nabla \eta\|_2\|\nabla (\eta ,u)\|_2 \|\nabla^3 u\|_0.  \label{202001052232}
\end{align}
Substituting the above estimate into \eqref{qwe201808100154628xxx}, using \eqref{aprpiosesnewxxxxx} and \eqref{aprpiosesnewxxxxzx},
we get \eqref{qwessebdaiseqinM0846dfgssgsd}.
\hfill$\Box$
\end{pf}
%%%%%%%%%%%%%%%%%%%%%%%%%%%%%%%%%%%%%%%%%%%%%%%%%%%%%%%
\begin{lem}\label{qwe201612132242nxsfssdfsxx}
  Under the assumptions \eqref{aprpiosesnewxxxxx}--\eqref{aprpiosesnewxxxxzx} with sufficiently small $\delta$, it holds that
\begin{align}
 & \frac{\mm{d}}{\mm{d}t}\left(  \bar{\rho} \sum_{|\alpha|=2}  \int \partial^\alpha \eta \partial^\alpha u\mm{d}y
 +\frac{\mu}{2}\|\nabla^3 \eta\|_0^2+\frac{\lambda}{2}\| \nabla^2 \mm{div}\eta \|_0^2\right)\nonumber \\
 & + P'(\bar{\rho})\bar{\rho}\|\nabla^2\mm{div}\eta\|^2_0   + {\kappa}  \|\nabla^3 \eta \|_0^2 /2\leqslant
    \bar{\rho} \|\nabla^2 u \|_{ 0}^2+ \delta c_{\mm{P}} \|\nabla^3 u\|_{0}^2.    \label{qweLem:03xx}
 \end{align}
\end{lem}
\begin{pf}
Applying $\partial^{\alpha}$ to \eqref{qwes0106pnnnn} and multiplying the resulting identity by $\partial^{\alpha}\eta$ in $L^2$,
one gets
\begin{align}
 &\frac{\mm{d}}{\mm{d}t}  \left(  \bar{\rho} \int \partial^\alpha \eta \partial^\alpha u\mm{d}y
 +\frac{\mu}{2}\|\nabla\partial^\alpha\eta\|_0^2+\frac{\lambda}{2}\|\mm{div}\partial^\alpha\eta\|_0^2\right)\nonumber \\
 & + P'(\bar{\rho})\bar{\rho}\|\mm{div}\partial^\alpha\eta\|^2_0
  + {\kappa} \|\nabla \partial^\alpha \eta \|_0^2 =\bar{\rho}\|\partial^\alpha u\|_0^2+I_6, \label{qwe2010154628}
\end{align}
where $|\alpha|=2$ and
$$I_6:=\int  \partial^\alpha {\mathcal{N}}^2  \cdot   \partial^\alpha  \eta\mm{d}y.$$

Analogously to \eqref{202001052232}, the integral $I_6$ can be bounded as follows.
\begin{align}
|I_6|\leqslant \|\nabla{\mathcal{N}}^2\|_0\|\nabla^3\eta\|_{0}
\lesssim_{\mm{P}} \|\nabla \eta\|_2^2\|\nabla (\eta ,u)\|_2.      \nonumber
\end{align}
Inserting the above estimate into \eqref{qwe2010154628}, and using \eqref{aprpiosesnewxxxxx}, \eqref{aprpiosesnewxxxxzx} and
Young's inequality, we obtain \eqref{qweLem:03xx}. \hfill$\Box$
 \end{pf}

With the help of Lemmas \ref{qwe201612132242nn}--\ref{qwe201612132242nxsfssdfsxx}, we can use the same arguments as for \eqref{201811262026}
 to deduce that there are constants $\delta_1^{\mm{P}}\in (0,1]$, ${c}_4^{\mm{P}}$, ${c}_5^{\mm{P}}$ and ${c}_{\mm{P}}$,
 such that for any $\delta\in (0,\delta_1^{\mm{P}}]$,
\begin{align}
&\frac{\mm{d}}{\mm{d}t}  {\mathcal{E}}_2^{{P}}+ c_{\mm{P}} \|\nabla ( u,\sqrt{\kappa}\eta)\|_2^2\leqslant 0,  \label{202001061326}
\end{align}
where
$$
\begin{aligned}
{\mathcal{E}}_2^{{P}}:= & {c}_4^{{P}}(\bar{\rho} \| \nabla^2 u\|_0^2+P'(\bar{\rho})\bar{\rho}\|\nabla^2 \mm{div}\eta\|^2_0
+ \kappa \| \nabla^3 \eta \|_0^2) + \bar{\rho} \sum_{  |\alpha| = 2}\int \partial^\alpha \eta \partial^\alpha u\mm{d}y \\
& +\frac{\mu}{2}\|\nabla^3 \eta\|_0^2+\frac{\lambda}{2}\|\nabla^2\mm{div}\eta\|_0^2 + {c}_5^{{P}}\left(\bar{\rho}\|u\|_0^2
+ {\bar{\rho}}\int\int_{\bar{\rho}/4}^{\bar{\rho}J^{-1}}\frac{P(z)}{z^2}\mm{d}z\mm{d}y + \kappa\|\nabla\eta \|_0^2\right),
\end{aligned}
$$
and ${\mathcal{E}}_2^{\mm{P}}$ satisfies
\begin{equation}
\nonumber
 \|(  u,\sqrt{\kappa}\nabla \eta)\|_2^2   \lesssim_{\mm{P}}  {\mathcal{E}}_2^{\mm{P}} \lesssim_{\mm{P}}  \|(  u,\sqrt{\kappa}\nabla \eta)\|_2^2   .
\end{equation}

 An integration of \eqref{202001061326} over $(0,t)$ implies that under the assumptions
\eqref{aprpiosesnewxxxxx}--\eqref{aprpiosesnewxxxxzx} with $\delta\leqslant \delta_1^{\mm{P}}$, one has
\begin{equation}
\label{2018112sfad62026xx}
  \|(  u,\sqrt{\kappa}\nabla \eta)\|_2^2 +
 \int_0^t\| \nabla (u, \sqrt{\kappa}\eta) \|_2^2 \mm{d}\tau  \leqslant K_{\mm{P}}^2/4\quad\mbox{ for any }t\geq 0,
 \end{equation}
where   ${c}_1^{\mm{P}}\geqslant 1$ and
\begin{equation}
\label{201saf911262060x}
K_{\mm{P}}:= 2\sqrt{    {c}_1^{\mm{P}}I_0^{\mm{h}}(u^0, \eta^0)}  .
\end{equation}

Similarly to Proposition \ref{pro:0401nxdxx}, we have the following local well-posedness result for the initial value problem \eqref{n0101nnn}--\eqref{2201912032006}:
\begin{pro} \label{qwepro:0401nxdxx}
Let $B_1^{\mm{P}}>0$ be a given constant, $( \eta^0,u^0)\in H^3 \times H^2$ and
$\|u^0\|_2\leqslant B_1^{\mm{P}}$. There is a sufficiently small $\delta_2^{\mm{P}}\in (0, 1]$,
such that if $\eta^0$ satisfy $\|\nabla \eta^0\|_2^2\leqslant \delta_2^{\mm{P}}$, then $3/4\leqslant\det(\nabla \eta^0+I)\leqslant 5/4$,
and there exists a unique strong solution $( \eta, u)\in C^0([0,T ), H^3 )\times\mathfrak{U}_T$ to the initial value problem  \eqref{n0101nnn}--\eqref{2201912032006} for some $T >0$, depending possibly on $B_1^{\mm{P}}$, $P(\cdot)$, $\bar{\rho}$, $\mu$
and $\lambda$. Moreover, $1/ 2\leqslant J\leqslant 3/2$.
\end{pro}
\begin{pf} Proposition \ref{qwepro:0401nxdxx} can be easily obtained by employing the same arguments as for
Proposition \ref{pro:0401nxdxx}, and hence its proof will be omitted here.    \hfill $\Box$
\end{pf}

Finally, with the help of \eqref{2018112sfad62026xx} and Proposition \ref{qwepro:0401nxdxx}, we can argue in a manner similar to
 the proof of Theorem \ref{201904301948} to easily get Theorem \ref{x2019saf04301948}.

 Next, we turn to the proof of Theorem \ref{20191asfas2041028}. An integration of \eqref{n0101nnn}$_2$ over $\mathbb{T}^3$ results in
 $\frac{\mm{d}}{\mm{d}t} \int u\mm{d}y=0$.
Thus, employing arguments similar to those in Section \ref{202001192326}, we can deduce that $(\eta,u)$ satisfies
\eqref{201905041x053}--\eqref{20191211197} with ${c}_3$ replaced by some positive constant ${c}_3^\mm{P}$.

We now prove \eqref{201905safa041053xxx}. Let $(\eta^0,u^0)$ be given in Theorem \ref{x2019saf04301948}, then there is
a unique global solution $(\eta^{\mm{l}},u^{\mm{l}} )\in  C^0(\mathbb{R}_0^+, H^{3})\times \mathfrak{U}_\infty $ to
\eqref{01dsaf16asdfsadffasfxx} with initial data $(\eta^0,{u}^0)$, satisfying
%% Moreover, $(\eta^{\mm{l}},u^{\mm{l}})$ enjoys the estimate
\begin{align}
  \|(  \bar{u}^{\mm{l}},\sqrt{\kappa}\nabla \eta^{\mm{l}})\|_2^2 +\int_0^t( \| \bar{u}^{\mm{l}}\|_3^2+{\kappa}\|\nabla \eta^{\mm{l}}\|_2^2 )e^{ {c}_3^\mm{P}(\tau-t)}\mm{d}\tau  \lesssim  e^{- {c}_3^\mm{P} t}   I_0^{\mm{h}}(  \bar{u}^0, {\eta}^0),
\label{201811sdf26asdfaf2026}
\end{align}
where $\bar{u}^{\mm{l}}:=u^{\mm{l}}-( {u}^0)_{\mathbb{T}^3}$.

Subtracting the initial value problem \eqref{01dsaf16asdfsadffasfxx} from \eqref{n0101nnn}--\eqref{2201912032006}, one gets
 \begin{equation}\label{01dsaf16sasadfasffafasdfasfxx}
                              \begin{cases}
\eta_t^{\mm{d}}=u^{\mm{d}}, \\[1mm]
 \bar{\rho} u_t^{\mm{d}}-\nabla( P'(\bar{\rho})\bar{\rho}\,\mm{div}\eta^{\mm{d}}+ \lambda\mm{div}u^{\mm{d}}  ) -\Delta ( \mu u^{\mm{d}}+
{\kappa} \eta^{\mm{d}})= \mathcal{N}^2, \\[1mm]
(\eta^{\mm{d}},u^{\mm{d}})|_{t=0}=0.
\end{cases}
\end{equation}
So, we employ arguments similar to those in deriving Lemmas \ref{qwe201612132242nn}--\ref{qwe201612132242nxsfssdfsxx} to arrive at
\begin{align}
  & \frac{\mm{d}}{\mm{d}t}  \left(  \bar{\rho} \|\nabla^i u^{\mm{d}}\|_0^2+ P'(\bar{\rho})\bar{\rho}\| \nabla^i\mm{div}\eta^{\mm{d}}\|^2
  + \kappa \| \nabla^{i+1}\eta^{\mm{d}} \|_0^2\right)+
 \mu \|\nabla^{i+1} u^{\mm{d}} \|_0^2   + \lambda \|\nabla^{i}\mm{div} u^{\mm{d}} \|_0^2 \nonumber \\
 & \quad \lesssim_{\mm{P}}  \|\nabla \eta\|_2\|\nabla^i u^{\mm{d}}\|_1 \|\nabla (\eta,u)\|_2, \quad i=0\mbox{ and }2,
\label{2019122asdfsaf51228}\\[2mm]
&  \frac{\mm{d}}{\mm{d}t}\left(  \bar{\rho} \sum_{|\alpha|=2} \int \partial^\alpha \eta^{\mm{d}} \partial^\alpha u^{\mm{d}}\mm{d}y
   + \frac{\mu}{2}\|\nabla^3 \eta^{\mm{d}}\|_0^2+\frac{\lambda}{2}\| \nabla^2 \mm{div}\eta^{\mm{d}} \|_0^2\right) \nonumber \\
& \qquad +P'(\bar{\rho})\bar{\rho}\|\nabla^2\mm{div}\eta\|^2_0 +  {\kappa}  \|\nabla^3\eta^{\mm{d}} \|_0^2/2
 \leqslant \bar{\rho}\|\nabla^2 u\|_2^2+c_{\mm{P}}  \|\nabla \eta\|_2 \|\nabla \eta^{\mm{d}}\|_2\|\nabla(\eta, u)\|_2.
 \label{Lem:031sfdaxx}
 \end{align}

 In addition, we can derive from \eqref{01dsaf16sasadfasffafasdfasfxx} that $(\eta^{\mm{d}})_{\mathbb{T}^3}=(u^{\mm{d}})_{\mathbb{T}^3}=0$.
Thus, similarly to \eqref{201811sdf262026},  we further get from \eqref{2019122asdfsaf51228}--\eqref{Lem:031sfdaxx} that
\begin{align}
&\frac{\mm{d}}{\mm{d}t} {\mathcal{E}}_2^{\mm{P,d}}+  {c}_3^\mm{P}({\mathcal{E}}_2^{\mm{P,d}}
+\|u^{\mm{d}} \|_3^2+ \kappa \| \nabla \eta^{\mm{d}}\|_2^2 ) \lesssim_{\mm{P}}
  \|\nabla \eta\|_2\|\nabla (\eta^{\mm{d}}, u^{\mm{d}})\|_2 \|\nabla (\eta,u)\|_2,
\label{201811sdsafdaff262026}
\end{align}
where ${\mathcal{E}}_2^{\mm{P,d}}$ is defined as ${\mathcal{E}}_2^{\mm{P}}$ with $(\eta ,u )$ replaced by $(\eta^{\mm{d}},u^{\mm{d}})$,
and satisfies the estimate
$$ \|  u^{\mm{d}}\|_2^2 +\kappa\| \eta^{\mm{d}}\|_3^2 \lesssim_{\mm{P}}  {\mathcal{E}}_2^{\mm{P,d}}
\lesssim_{\mm{P}} \|  u^{\mm{d}}\|_2^2 +\kappa\| \eta^{\mm{d}}\|_3^2 . $$

Noting that $(\eta,u)$ also satisfies \eqref{201905041x053} with $ {c}_3^\mm{P}$ in place of $ {c}_3$,
integrating \eqref{201811sdsafdaff262026} over $(0,t)$ and then using \eqref{201905041x053} with ${c}_3$ replaced by ${c}_3^\mm{P}$,
and \eqref{201811sdf26asdfaf2026}, we immediately obtain the desired estimate \eqref{201905safa041053xxx}.
 This completes the proof of Theorem \ref{20191asfas2041028}.

\section{Proof of Proposition \ref{pro:0401nxdxx}}\label{201912062021}

In this section we give the proof of Proposition \ref{pro:0401nxdxx} which is divided into four subsections.

\subsection{Solutions of the linearized problem}\label{2020010813027}

   For given $(\eta^0,w)$, consider the solution of the initial value problem:
\begin{equation}
\label{201912060857}
  \begin{cases}
 \eta_t=u  ,\\ {\rho} u_t+\nabla_{\ml{A}}q-\mu \Delta_{\ml{A}}u=
{\kappa}\Delta \eta  ,    \\
\mm{div}_{\mathcal{A}}u=0 ,\\
u|_{t=0}=u^0 ,
  \end{cases}
  \end{equation}
where $\ml{A}$ is defined by $\zeta$ and
\begin{align}
\zeta=\int_0^t w\mm{d}y+\eta^0+y.
\label{202001201942}
\end{align}
Then we have the following conclusion:
\begin{pro} \label{qwepro:0sadfa401nxdxx}
Let $B_2\geqslant 1$, $(\eta^0,u^0)\in H^3\times H^2$. Assume
\begin{align}
&w\in  C^0(\bar{I}_T,H^2 )\cap L^2(I_T,H^3),\ w_t\in L^2(I_T,H^1),\ w|_{t=0}=u^0,\label{201912282123} \\
& \sqrt{\|\nabla w\|_{C^0(\overline{I_T},H^1)}^2 +\|\nabla w\|_{L^2(I_T,H^2)}^2
+\|\nabla w_t\|_{L^2(I_T,L^2)}^2}\leqslant B_2,\label{2019122821231}
\end{align}
then there are a sufficiently small constant $\delta_3\in (0,1]$ independent of any parameters, and a local existence time
\begin{align}
T\in (0,\min\{1,(\delta_3/2B_2)^4\}],  \label{202001212374}
\end{align}
such that for any $\eta^0$ satisfying $\|\nabla \eta^0\|_2 \leqslant \delta_3$, there exists a unique local strong solution
 $(u,q)\in\mathfrak{U}_T\times L^\infty(I_T, \underline{H}^2£©$ to \eqref{201912060857}. Moreover,
\begin{align}
& \|\tilde{\mathcal{A}}\|_2\lesssim_0 \|\nabla \eta\|_2\lesssim_0 \delta,\ \|\tilde{\mathcal{A}}_t\|_1\lesssim_0 \| \nabla w\|_1, \nonumber \\
& \|    u  \|_{\mathfrak{U}_T} + \|{q}\|_{\mathfrak{Q}_T}  \leqslant  c_\kappa(  1+\|u^0\|_2^3)+  \|\nabla w\|_{C^0(\overline{I_T},H^1)}/2  ,\label{201912242055}\\
&  \|\nabla q \|_{L^\infty({I}_T,H^1)}\lesssim_\kappa B_2  \sqrt{1+\|u^0\|_2^4+ B_2(1+ \| u^0\|_2^3) } +1 . \label{201912242055safdsaf}
\end{align}
\end{pro}
\begin{pf}
 Let $\eta:=\zeta-y$ with $\zeta$ being given by \eqref{202001201942}, $\|\nabla \eta^0\|_2\leqslant \delta/2$ and $T$ satisfy
\begin{align}
\label{201912061028}
T\in (0,1]\;\;\mbox{ and }\;\; 2{T}^{1/4} B_2 \leqslant \delta \in (0,1].
\end{align}
By virtue of \eqref{2019122821231} and \eqref{201912061028}, it is easy to check that $\nabla \eta\in C^0(\bar{I}_T,H^2 )$ and
\begin{equation}
\label{201912061030}
\|\nabla \eta(t)\|_2 \leqslant\delta/2+\sqrt{t} \|\nabla w\|_{L^2(I_T,H^2)}\leqslant\delta\;\;\mbox{ for any }t\in I_T.
\end{equation}
Moreover, by \eqref{201912061030} we have for sufficiently small $\delta$ that
\begin{align}
 & 1/2 \leqslant J\leqslant 3/2, \label{2019123023dasa17} \\
& \|\tilde{\mathcal{A}}_t\|_j\lesssim_0 \| \nabla w\|_j\mbox{  for }1\leqslant j\leqslant 2, \label{201912221512}\\
& \|\tilde{\mathcal{A}}\|_2 \lesssim_0   \|\nabla \eta\|_2. \label{201912301350}
\end{align}

Inspired by the proof of \cite[Theorem 4.3]{GYTILW1}, our next step is to solve the linear problem \eqref{201912060857}
by applying the Galerkin method.
The space $H^2_\sigma$ possesses a countable orthogonal basis $\{\varphi^i\}_{i=1}^\infty$, since it is a separable Hilbert space.
For each $m\geqslant 1$ we define $\psi^i=\psi^i(t):=\nabla\zeta\varphi^i$, where $\zeta=\eta +y$. Then, by \eqref{201912061030},
$\psi^i(t)\in \mathfrak{H}(t):= \{u\in H^2~|~\mm{div}_{\mathcal{A}}u=0\}$
and $\{\psi^i(t)\}_{i=1}^\infty$ is a basis of $\mathfrak{H}(t)$ for each $t\in \overline{I_T}$. Moreover,
\begin{equation}
\label{201912102051}
\psi_t^i =R \psi^i,
\end{equation}
where $R:=\nabla w\mathcal{A}^{\mm{T}}$.
Thanks to \eqref{201912061030}, \eqref{201912221512} and \eqref{201912301350}, we have upper bounds:
\begin{align}
&\|R\|_j\lesssim_0  \|\nabla w\|_j\;\;\mbox{ for }0\leqslant j\leqslant 2,\label{20asfd1912241400} \\
& \|R\|_{L^\infty}\lesssim_0 \|\nabla w\|_{L^\infty}\lesssim_0 \|\nabla w\|_1^{1/2}\|\nabla w\|_2^{1/2},\label{201912241400} \\
& \|\nabla  R \|_{L^3}\lesssim_0 \|\nabla R\|_0^{1/2}\|\nabla R\|_{L^6}^{1/2}\lesssim_0
 \|\nabla w\|_{1}^{1/2}\|\nabla w\|_{2}^{1/2} , \label{201912301544}\\
 & \|R_t\|_0\lesssim \|\nabla w_t\|_0+\|\nabla w\|_1^2.\label{201912301544xx}
 \end{align}

For any integer $m\geqslant 1$, we define the finite-dimensional space
$ \mathfrak{H}_m(t):=\mm{span}\{\psi^1,\ldots, \psi ^m\}\subset \mathfrak{H}(t)$,
and write $\mathcal{P}^m(t):\mathfrak{H}(t)\to \mathfrak{H}_m(t) $ for the $\mathfrak{H}$ orthogonal projection
onto $ \mathfrak{H}_m(t)$. Clearly, for each $v\in \mathfrak{H}(t)$,
$\mathcal{P}_m(t) v \rightarrow v $ as $m\rightarrow\infty$ and $\|\mathcal{P}_m(t) v \|_2\leqslant \|v\|_2$.

Now we define an approximate solution
$$u^m(t)=a_j^m(t) \psi^j\;\;\mbox{ with }a_j^m:\ \overline{I_T}\rightarrow  \mathbb{R}\mbox{ for }j=1,\ldots, m,$$
where $m\geqslant 1$ is given.
We want to choose the coefficients $a_j^m$, so that for any $1\leqslant i\leqslant m$,
\begin{equation}\label{appweaksolux}
 \rho \int  u ^m_t\cdot \psi^i \mm{d} y +\mu\int
\nabla_{\mathcal{A}}  u ^m:\nabla_{\mathcal{A}} \psi^i \mm{d} y = \kappa \int \Delta \eta  \cdot \psi^i \mm{d} y
\end{equation}
with initial data
$$ u ^m(0)=\mathcal{P}^m u _0\in \mathfrak{H}_m.$$

Let
$$  \begin{aligned}
&A:=(a_1^m, ,\ldots,  a_m^m)^{\mm{T}},\quad X:=\left( \rho \int\psi^i\cdot \psi^j \mm{d} y\right)_{m\times m},\\
 &Y:= \kappa \left(\int \Delta \eta \cdot\psi^1 \mm{d} y,\ldots,  \int \Delta \eta \cdot\psi^m \mm{d} y\right)^{\mm{T}},\\
&Z:=\left( \int \left(\rho R \psi^i\cdot \psi^j+
\mu\nabla_{\mathcal{A}} \psi^i:\nabla_{\mathcal{A}} \psi^j \right)\mm{d} y\right)_{m\times m}.
\end{aligned}$$
By the regularity of $w$ in \eqref{201912282123}, we have
\begin{align}
\label{202001200222}
X,\ Y\in C^{1,1/2}(\overline{I_T}),\ Z\in C^{0,1/2}(\overline{I_T})\mbox{ and }Z_t\in L^2(I_T).
\end{align}

 Noting that $X$ is invertible, we can thus make use of \eqref{201912102051} to rewrite \eqref{appweaksolux} as follows.
 \begin{equation}  \label{appweaksolu}
  A _t = X^{-1}(Y-ZA) .
\end{equation}
By virtue of the well-posedness theory of ODEs (see \cite[Section 6 in Chapter II]{WWODE148}), the equation \eqref{appweaksolu} has
a unique solution $A\in C^{1,1/2}(\overline{I_T})$. Thus, we have gotten the existence of the approximate solution $u^m(t)=a_j^m(t)\psi^j$.
Moreover, it is easy to see that $\ddot{a}_j^m(t)\in L^2({I_T})$ by the regularity \eqref{202001200222}.
Next, we derive uniform-in-$m$ estimates for $u^m$.

Thanks to \eqref{20160614fdsa1957}, we easily derive from \eqref{appweaksolux} with $u^m$ in place of $\psi$ that
for sufficiently small $\delta$,
\begin{equation}
\label{201912112038}
\frac{\mm{d}}{\mm{d}t}\| u ^m \|_0^2 + c_0\|\nabla u^m\|_0^2 \lesssim_0 \|\nabla \eta\|_0^2 .
\end{equation}

By \eqref{201912102051},
\begin{equation}
\label{201912122222}
u^m_t- R  u^m= \dot{a}_i^m  \psi^i.
\end{equation}
Then we can replace $\psi$ by $ \dot{a}_i^m  \psi^i$ in \eqref{appweaksolux}  to deduce that
\begin{align}
& \rho\|u ^m_t\|_0^2+\mu\int
\nabla_{\mathcal{A}}  u ^m:\nabla_{\mathcal{A}} u ^m_t \mm{d}  =\rho \int u_t^m\cdot(R u^m)\mm{d}y y\nonumber  \\
& \quad +\mu \int \nabla_{\mathcal{{A}}}u^m:\nabla_{\mathcal{A}} (R u^m)\mm{d}y
 +\kappa\int \Delta \eta  \cdot \left(   u ^m_t- R u^m\right)\mm{d}y. \label{appweakssafolux}
\end{align}
From \eqref{appweakssafolux} we can further deduce that
\begin{align}
& {\mu } \frac{\mm{d}}{\mm{d}t} \|\nabla_{\mathcal{A}}  u ^m\|_0^2
+  \rho\|u ^m_t\|_0^2 \lesssim_\kappa  \|(R u^m,\Delta \eta)\|_{0}^2+\|\nabla u^m\|_0(\|\nabla(R u^m)\|_0
+\|\mathcal{A}_t\|_{L^\infty} \|\nabla u^m\|_0 )\nonumber \\
& \qquad  \lesssim_\kappa \|\nabla w\|_2(1+\|\nabla w\|_1 )   (\|u^m\|_0^2+\|\nabla_{\mathcal{A}}u^m\|_0^2)+\|\nabla \eta\|_1^2,
\label{appwfolux}
\end{align}
where we have used  \eqref{201905041149}, \eqref{201912221512},  \eqref{201912241400} and \eqref{201912301544}   in the last inequality.

Summing up \eqref{201912112038} and \eqref{appwfolux}, and using Gronwall's lemma, \eqref{201912061028} and \eqref{201912061030},
we infer that for any $t\in\overline{I_T}$,
\begin{align}
&\| u ^m
\|_1^2   + \int_0^t\|(\nabla u^m,u ^m_t)\|_0^2\mm{d}\tau \lesssim_\kappa (1+\|\mathcal{P}^mu^0\|_1^2)
 e^{\int_0^t \|\nabla w\|_2(1+\|\nabla w\|_1 )  \mm{d}\tau} \nonumber \\
&\qquad \quad \times\left(1+ \int_0^t  \|\nabla w\|_2(1+\|\nabla w\|_1 ) \mm{d}\tau\right)\lesssim  1+\|u^0\|_2^2.
\label{201912241359}
\end{align}

In view of \eqref{201912102051}, we get from \eqref{appweaksolux} that
 \begin{align}
&  \rho \int  u ^m_{tt}\cdot \psi \mm{d} y +\mu\int
\nabla_{\mathcal{A}}  u ^m_t:\nabla_{\mathcal{A}} \psi \mm{d} y  % \nonumber \\
=  \kappa \int (\Delta w\cdot\psi + \Delta \eta\cdot (R \psi))\mm{d} y
-\rho \int  u ^m_{t}\cdot (R \psi) \mm{d} y \nonumber \\
&\quad -\mu\int
(\nabla_{\mathcal{A}_t} u^m:\nabla_{\mathcal{A}} \psi
+\nabla_{\mathcal{A}} u^m:(\nabla_{\mathcal{A}_t} \psi+ \nabla_{\mathcal{A}}(R\psi ))\mm{d} y\;\;\mbox{  a.e. in }I_T. \label{201912102242}
\end{align}

Noting that (referring to \cite[Theorem 1.67]{NASII04})
$$
\begin{aligned}
& \frac{1}{2}\|u ^m_t\|_0^2
 - \int  u_{t}^m\cdot (R u^m)\mm{d}y- \left.\left(\frac{1}{2}\|u ^m_t\|_0^2
 - \int  u_{t}^m\cdot ( R u^m)\mm{d}y\right)\right|_{t=0}\\
 & = \int_0^t \left(\int  u ^m_{\tau\tau}\cdot  ( u^m_{\tau}-Ru^m) \mm{d} y -\int  u ^m_{\tau}\cdot  (Ru^m)_{\tau} \mm{d} y\right)\mm{d}\tau,
 \end{aligned} $$
and utilizing \eqref{201912122222} and the above identity, we deduce from \eqref{201912102242} with $\psi$ replaced
by $(u^m_t-Ru^m)$ that
 \begin{align}
 & \rho\left(\frac{1}{2 }\|u ^m_t\|_0^2  - \int   u_{t}^m\cdot (R u^m)\mm{d}y\right)
 +\mu\int_0^t \|\nabla_{\mathcal{A}}  u ^m_\tau\|_0^2\mm{d}\tau \nonumber \\
& \qquad = \left.\rho\left(\frac{1}{2}\|u ^m_t\|_0^2 - \int u_{t}^m\cdot ( Ru^m)\mm{d}y\right)\right|_{t=0}+I_7,
\label{201912241249}
\end{align}
where
$$ \begin{aligned}
  I_7 :=&\int_0^t\left(\kappa\int (\Delta w\cdot (u^m_\tau-R u^m) +\Delta \eta\cdot (R (u_\tau^m- R u^m)))\mm{d}y \right.
  \nonumber \\
& \quad +\rho \int (u_{\tau}^m\cdot (R^2 u^m) - 2u_{\tau}^m\cdot (R u^m_{\tau}) - u_{\tau}^m\cdot (R_{\tau} u^m) )\mm{d}y \nonumber \\
&  -\mu\int (  \nabla_{\mathcal{A}} u^m: (\nabla_{\mathcal{A}_{\tau}} (u_{\tau}^m-R u^m)
+\nabla_{\mathcal{A}}(R (u_{\tau}^m -R u^m)))\nonumber \\
&\quad  +\nabla_{\mathcal{A}_{\tau}}  u^m :\nabla_{\mathcal{A}} (u_{\tau}^m-R u^m)
- \nabla_{\mathcal{A}}  u^m_\tau :\nabla_{\mathcal{A}} (R u^m)) \mm{d} y\bigg)\mm{d}\tau .
\end{aligned}   $$

Recalling that
\begin{align}
  \|\nabla w\|_0\leqslant \int_0^t\|\nabla w_{\tau}\|_0\mm{d}\tau+\|\nabla u^0\|_0, \nonumber
\end{align}
we make use of \eqref{201912061028}, and \eqref{20asfd1912241400}, \eqref{201912241359} and the above estimate to deduce
from \eqref{201912241249} that
 \begin{align}
 \|u ^m_t\|_0^2  +\int_0^t\|\nabla  u ^m_{\tau}\|_0^2\mm{d}\tau
& \lesssim_\kappa\|\nabla w\|_{0}\|\nabla w\|_1\|u^m\|_1^2+\|\nabla w^0\|_{1}^2\|u^0\|_2^2 +\|u^m_t|_{t=0}\|_0^2+ I_7 \nonumber \\
& \lesssim_\kappa   \|\nabla w\|_1(1+ \| u^0\|_2^3) + \| u^0\|_2^4 +\|u ^m_t|_{t=0}\|_0^2+ I_7,
\label{201912242022}
\end{align}
where the last two terms on the right hand of \eqref{201912242022} can be bounded as follows. Replacing $\psi$ by $(u_t^m-Ru^m)$ in \eqref{appweaksolux}, one sees that
\begin{align}
 \rho\| u^m_t \|_0^2 = & \kappa \int \Delta \eta  \cdot (u ^m_t-R u^m) \mm{d} y
 +\mu\int \Delta_{\mathcal{A}}u^m:(u_t^m-R u^m)\mm{d}y+\rho\int u_t^m\cdot (R u^m)\mm{d}y, \nonumber
\end{align}
whence,
$$ \|u^m_t \|_0^2\lesssim_\kappa\|\nabla (\eta,u^m)\|_1^2 +\|\nabla w\|_1^2\|u^m\|_1^2\;\;\;\mbox{ for any }t\in [0,T). $$
In particular, \begin{equation}
\label{201912112039}
\| u ^m_t|_{t=0} \|_0^2\lesssim_\kappa \|\nabla \eta^0\|_1^2+(1+\|u^0\|_2^2)^2 \lesssim_\kappa 1+\|u^0\|_2^4 .
\end{equation}

%Let $e_j$ be the unit vetoer with the $j$-th component being $1$ for $1\leqslant j\leqslant 3$. Since $e_i\in \mathfrak{H}(t)$ can be expressed by  $ \{ \psi^i \}_{i=1}^\infty$, thuswe can replace $\psi$ by $e_j$ in  \eqref{appweaksolux}, and get $$\int  u ^m_t   \mm{d} y=0 .$$By Poinc\'are's inequality, we have $$\|u_t^m\|_{L^6}\lesssim \|u_t^m\|_{1}\lesssim  \|\nabla u_t^m\|_{0}.$$
On the other hand, the last term on the right hand of \eqref{201912241249} can be estimated as follows.
\begin{align}
 I_7\lesssim_\kappa &  \int_0^t(
 \|u^m\|_1(\|R \|_1\|\nabla w\|_1+\|R \|_1^2 +
 \|R_\tau \|_0 (\|u_\tau^m\|_{0}+\|u_\tau^m\|_{0}^{1/2}\|\nabla u_\tau^m\|_{0}^{1/2}) + (\|R\|_{L^\infty}
\nonumber \\
&\quad + \|\nabla R\|_{L^3
 }+\|  {\mathcal{A}_\tau}\|_{1}^{1/2}\|  {\mathcal{A}_\tau}\|_{2}^{1/2})(\| \nabla u_\tau^m\|_0
 +\|  u^m\|_1(\|R\|_{L^\infty}+\|\nabla R\|_{L^3}  )))\nonumber \\
 &\quad  +\| u^m_\tau \|_0(  \| \nabla w\|_1+\|R \|_1 +(\| R \|_{1}^2 +\|\nabla R\|_{L^3}) \|u^m\|_1) +
\|u_\tau^m\|_0^2 \|R\|_{L^\infty} )\mm{d}\tau \nonumber \\
 \lesssim_\kappa  &(1+\|u^0\|_2^2)\left(1+ \sup_{t\in \overline{I_T}}\|u_\tau^m\|_0+\left(\int_0^t\| \nabla u_\tau^m\|_0^2\mm{d}\tau\right)^{1/2}  \right)+ \delta \sup_{t\in \overline{I_T}}\|u_\tau^m\|_0^2,  \label{202001221640}
\end{align}
where one has used \eqref{201912061028}, \eqref{201912221512} and \eqref{20asfd1912241400}--\eqref{201912301544xx}
in the second inequality in \eqref{202001221640}.

Putting \eqref{201912112039} and \eqref{202001221640} into \eqref{201912242022}, we get for sufficiently small $\delta$ that
\begin{align}
 \sup_{ t\in \overline{I_T}}\| u ^m_t \|_0^2 +   \| u ^m_{t}\|_{L^2(I_T,H^1)}^2 \lesssim_\kappa  1+\|u^0\|_2^4
 + \|\nabla w\|_{C^0(\overline{I_T},H^1)}(1+ \| u^0\|_2^3) .\label{appwasdfasfsadfeaksolux}
\end{align}
Summing up \eqref{201912241359} and \eqref{appwasdfasfsadfeaksolux}, we conclude
\begin{align}
 \sup_{ t\in \overline{I_T}}( \| u ^m \|_1^2+\|u ^m_t \|_{0}^2)+ \|(u ^m,u _t^m)\|_{L^2(I_T,H^1)}^2
\lesssim_\kappa  1+\|u^0\|_2^4+  \|\nabla w\|_{C^0(\overline{I_T},H^1)}(1+ \| u^0\|_2^3). \label{sumestimeat}
\end{align}

In view of \eqref{sumestimeat}, up to the extraction of a subsequence, we have, as $m\to\infty$,
$$\begin{aligned}
& (u^m,u^m_t)\rightarrow (u,u_t) \mbox{ weakly-* in }L^\infty(I_T,H^1\times L^2 ) ,\\
& (u^m,u^m_t)\rightarrow (u,u_t) \mbox{ weakly  in }L^2(I_T,H^1\times H^1),\\
& u (0)= u _0 \;\; \mbox{ and }\; \mm{div}_{\mathcal{A}} u =0.
\end{aligned}$$
Moreover,
\begin{align}
&  \| u_t \|_{L^\infty(I_T,L^2)}^2+  \|u_t\|_{L^2(I_T,H^1)}^2  \lesssim_\kappa 1+\|u^0\|_2^4
 +  \|\nabla w\|_{C^0(\overline{I_T},H^1)}(1+ \| u^0\|_2^3).  \label{2019121asdfa12157}
   \end{align}
Thus, taking to the limit as $m\to\infty$ in \eqref{appweaksolux}, we obtain that for any $\phi \in \mathfrak{H}$,
\begin{equation}\label{0507n1}
\rho \int u _t\cdot \phi \mm{d} y +\mu\int
\nabla_{\mathcal{A}} u:\nabla_{\mathcal{A}} \phi\mm{d} y=\kappa\int\Delta\eta\cdot\phi\mm{d} y\;\; \mbox{ a.e. in }I_T,
\end{equation}which can be regarded as a weak solution of the following Stokes equations
 \begin{equation*}
 \nabla_{\mathcal{A}} {q} -\mu \Delta_{\mathcal{A}} u = {F}:= \kappa\Delta\eta -\rho u_t\in L^\infty(I_T,L^2)\cap L^2(I_T,H^1) . \end{equation*}

Next, we want to show spatial regularity of $u$. For this purpose, we further assume that $\delta $ is so small that
$\eta$ satisfies \eqref{2312018031adsadfa21601xx}, we also refer the reader to the proof of (4.3) in \cite[Lemama 4.2]{JFJSOMITNx}.

Let $\tilde{F}=F(\zeta^{-1},t)$, then $\tilde{F}$ is also a periodic function defined on $\mathbb{T}^3$. Moreover,
by virtue of \eqref{2019123023dasa17},
$$  \mm{ess} \sup_{ t\in {I_T}}\| F \|_{0}^2+  \int_0^T\|  F\|_{1}^2\mm{d}\tau<\infty. $$
 Applying the regularity theory of the Stokes problem, we have that for any given $t\in I_T$, there is a unique strong solution
 $v(x,t)\in H^3$ with a unique associated function $p(y,t)\in\underline{H}^3$, such that
 \begin{equation}\label{appstokes}
  \begin{cases}
\nabla  p -\mu \Delta v =\tilde{F}  , \\
\mm{div}v =0.
  \end{cases}
 \end{equation}
 Moreover, for $k=0$ and $1$ it holds that $\|v\|_{2+k} + \|\nabla p\|_{k}\lesssim \| \tilde{F}\|_{k}$.

Let $w=v(\zeta,t)$ and $q=p(\zeta,t)-( p(\zeta,t))_{\mathbb{T}^3}$, then $(w,q)\in H^3\times\underline{H}^2$ for given $t\in I_T$.
Moreover, by \eqref{appstokes} one sees
 \begin{equation}
 \label{appsdfstokes}
 \begin{cases}
\nabla_{\mathcal{A}}  {q} -\mu \Delta_{\mathcal{A}} w= {F}  , \\
\mm{div}_{\mathcal{A}} w =0 .
  \end{cases}
 \end{equation}
Obviously, $u=w$ by \eqref{0507n1}.

 Now, we further rewrite \eqref{appsdfstokes} with $u$ in place $w$ as follows.
 \begin{equation}\nonumber
 \begin{cases}
\nabla  {q} -\mu \Delta {u}= G:={F} +
\mu(\mm{div}_{\tilde{\mathcal{A}}}\nabla_{\mathcal{A}} {u}+\mm{div} \nabla_{\tilde{\mathcal{A}}} {u} ) -\nabla_{\tilde{\mathcal{A}}}  {q} , \\
\mm{div} {u} =-\mm{div}_{\tilde{\mathcal{A}}} {u} .
  \end{cases}
 \end{equation}
Then, we apply the regularity theory for the Stokes equation to get that for sufficiently small $\delta$,
\begin{equation}  \|  u \|_{2+k}  +
 \| \nabla{q} \|_{k} \lesssim_\kappa \| u_t\|_{k} +\|\nabla \eta \|_{1+k} \lesssim_\kappa 1+\| u_t\|_{k}\mbox{ for }k=0,\ 1. \label{201912221611}
 \end{equation}  Consequently, combining \eqref{2019121asdfa12157} with \eqref{201912221611}, we arrive at the estimate
\begin{align}
 \|u \|_{\mathfrak{U}_T} + \|{q}\|_{\mathfrak{Q}_T}  \lesssim_\kappa 1+\|u^0\|_2^2+
\sqrt{\|\nabla w\|_{C^0(\overline{I_T},H^1)}(1+ \| u^0\|_2^3)}, \nonumber
\end{align}
which, together with Young's inequality, yields \eqref{201912242055}.
Moreover, analogously to \eqref{2017020614181721}, we obtain \eqref{201912242055safdsaf}.

Since $ u\in L^2(I_T, H^3)$ and $ u_t\in L^2(I_T, H^1)$, $u\in C^0(I_T, H^2)$ by means of a classical regularization method
(please refer to \cite[Lemma A.4]{GYTILW1}).

%Recalling \eqref{appstokes}$_1$ and the relation $ v_t=(u_t-u\cdot\nabla_{\mathcal{A}}u)|_{y=\zeta^{-1}}=u_t|_{y=\zeta^{-1}}-v\cdot\nabla v$, we have$$ \rho( v_t+v\cdot\nabla v)+\nabla  p -\mu \Delta v ={\kappa}\Delta \eta|_{y=\zeta^{-1}} .$$ Multiplying the above identity by $\phi\in H_\sigma^1$, and then differentiate the resulting identity, we get that\begin{equation}\label{absolution}\begin{aligned}\rho \frac{\mathrm{d}}{\mathrm{d}t}\int  v _t\cdot \phi \mm{d} x =\kappa\int \partial_t(\Delta \eta|_{y=\zeta^{-1}} )\cdot \phi\mm{d}x -\mu \int  (\nabla   v_t  :\nabla  \phi ) \mm{d} x-  \rho\int  (v\cdot\nabla v)_t  \cdot  \phi  \mm{d}x \end{aligned} \end{equation}   a.e. $t\in I_T$. Noting that $(\nabla v_t,\partial_t(\Delta \eta|_{y=\zeta^{-1}} ))\in L^\infty(I_T, L^2)$ and $v\in L^\infty(I_T, H^2)$, we get $v _{tt}\in L^2(I_T, {H}^{-1}_\sigma)$ from \eqref{absolution}. In addition, $v_t\in  L^2(I_T,  {H}_\sigma^1)$, thus $v _t\in C^0(\bar{I}_T,L^2)$, which implies that $u_t\in C^0(\bar{I}_T,L^r)$ for any $r\in [1,2)$. Moreover $q\in  C^0(\bar{I}_T,W^{1,p})$ by the relation \eqref{201912060857}$_2$.

Finally, the solution $(u,q)$ to the linear problem \eqref{201912060857} is obviously unique in the function class $\mathfrak{U}_T\times L^\infty(I_T,\underline{H}^2)$. This completes the proof of Proposition \ref{qwepro:0sadfa401nxdxx}. \hfill $\Box$
\end{pf}

\subsection{Proof of Proposition \ref{pro:0401nxdxx}}
Let $B_1>0$ be fixed, $(\eta^0,u^0)\in H^3 \times H^2 $, $\|u^0\|_2\leqslant B_1$ and $\|\nabla\eta^0\|_2^2\leqslant\delta\leqslant\delta_3$.
Denote
$$B_2:= 2c_\kappa(  1+B_1^3)+B_1,$$
 where the constant $c_\kappa $ comes from the estimate \eqref{201912242055}. By Proposition \ref{qwepro:0sadfa401nxdxx},
 we can construct a function sequence $\{ u^k,{q}^{k}\}_{k=1}^\infty$ defined on $Q_T$ with $T$ satisfying \eqref{202001212374}.
 Moreover,
\begin{itemize}
  \item for $k\geqslant 1$, $({u}^{k+1},q^{k+1})\in  \mathfrak{U}_T\times L^\infty(I_T ,\underline{H}^2) $ and
 \begin{equation}\label{iteratiequat}\begin{cases}
\rho u_t^{k+1}+\nabla_{\mathcal{A}^k} {q}^{k+1}-\mu \Delta_{\mathcal{A}^k} {u}^{k+1}=\kappa\Delta \eta^k ,\\[1mm]
\mathrm{div}_{\mathcal{A}^k} {u}^{k+1}=0\end{cases}  \end{equation}
with initial condition ${u}^{k+1} |_{t=0}=u^0$,
where $\eta^k:=\int_0^tu^k\mm{d}\tau+\eta^0$ and $\mathcal{A}^k$ is defined by $\zeta^k:=\eta^k+y$;
  \item $({u}^1,q^1)$ is constructed by Proposition \ref{qwepro:0sadfa401nxdxx} with $w=0$;
\item
the solution sequence $\{  {u}^k,q^k\}_{k=1}^\infty$ satisfies the following uniform estimates: For any $k\geqslant 1$,
\begin{align}
\label{n0531}
&\|\tilde{\mathcal{A}}^k\|_2\lesssim_0 \|\nabla \eta^k\|_2\lesssim_0 \delta,\ \| {\mathcal{A}}_t^k\|_1\lesssim_0 \| \nabla u^k\|_1 , \\  &\|u^k \|_{\mathfrak{U}_T}
  +\|{q}^k\|_{\mathfrak{Q}_T}  \leqslant   B_2,\ \|\nabla q^k \|_{L^\infty({I}_T,H^1)}\lesssim_\kappa 1+B_2^2\label{n053112}.\end{align}
\end{itemize}
In order to take limits in \eqref{iteratiequat}, we have to show $\{{u}^k,q^k\}_{k=1}^\infty$ is a Cauchy sequence.
To this end, we define for $k\geqslant 2$,
$$( \bar{ {u}}^{k+1},  \bar{q}^{k+1})= ({u}^{k+1}- {u}^k,{q}^{k+1} -q^{k}),$$
which satisfies
 \begin{equation}\label{difeequion} \begin{cases}
 \Delta  \bar{q}^{k+1} =\mathcal{M}_k  , \\[1mm]
\rho\bar{ { u}}_t^{k+1}+\nabla \bar{q}^{k+1}-\mu \Delta  \bar{ {u}}^{k+1}-\kappa \Delta \bar{\eta}^k=\mathcal{N}_k,\\[1mm]
\mathrm{div} \bar{u}^{k+1}= -(\mathrm{div}_{\tilde{\mathcal{A}}^k
-\tilde{\mathcal{A}}^{k-1}} {u}^{k+1}+\mathrm{div}_{\tilde{\mathcal{A}}^{k-1}} \bar{u}^{k+1})\\
\bar{ {u}}^{k+1} |_{t=0}=0, \end{cases}  \end{equation}
where
$$
\begin{aligned}
\mathcal{M}_k:=&{\kappa}\Delta(\mm{div}\bar{\eta}^k+ \mm{div}_{\tilde{\mathcal{A}}^{k}
-\tilde{\mathcal{A}}^{k-1}}  {\eta}^{k}+\mm{div}_{ \tilde{\mathcal{A}}^{k-1}}  \bar{\eta}^{k})\nonumber  \\
& +\rho (\mm{div}_{\mathcal{A}_t^k-\mathcal{A}_t^{k-1}}{u}^{k+1} + \mm{div}_{\mathcal{A}_t^{k-1}}\bar{u}^{k})
-(\mm{div}_{\tilde{\ml{A}}^k-\tilde{\ml{A}}^{k-1}}\nabla_{\ml{A}^k}{q}^{k+1} +
\mm{div}_{\tilde{\ml{A}}^{k-1}}\nabla_{\tilde{\ml{A}}^k-\tilde{\ml{A}}^{k-1}} {q}^{k+1}\nonumber &\\
&+ \mm{div}_{\tilde{\ml{A}}^{k-1}}\nabla_{\tilde{\ml{A}}^{k-1}}\bar{q}^{k+1}
+ \mm{div} (\nabla_{\tilde{\ml{A}}^{k}-\tilde{\ml{A}}^{k-1}}q^{k+1}+ \nabla_{\tilde{\ml{A}}^{k-1}}\bar{q}^{k+1})), \\
\mathcal{N}_k:=& \mu ( \mm{div}_{\tilde{\ml{A}}^k-\tilde{\ml{A}}^{k-1}}\nabla_{\ml{A}^k}u^{k+1}
+ \mm{div}_{ \tilde{\ml{A}}^{k-1}}\nabla_{\ml{A}^k-\ml{A}^{k-1}}u^{k+1}+\mm{div}_{ \tilde{\ml{A}}^{k-1}}\nabla_{\ml{A}^{k-1}}\bar{u}^{k+1}\\
&+ \mm{div}\nabla_{\tilde{\ml{A}}^k-\tilde{\ml{A}}^{k-1}}u^{k+1}+ \mm{div} \nabla_{\tilde{\ml{A}}^{k-1}}\bar{u}^{k+1})
-(\nabla_{\tilde{\ml{A}}^k-\tilde{\ml{A}}^{k-1}}q^{k+1}+\nabla_{\tilde{\ml{A}}^{k-1}}\bar{q}^{k+1}).
\end{aligned}$$

Keeping in mind that for sufficiently small $\delta$,
$$
\begin{aligned}
& \| {\ml{A}}^k- {\ml{A}}^{k-1}  \|_2\lesssim_0 \|\nabla \bar{\eta}^k\|_2\lesssim_0 T^{1/2} \| \nabla \bar{u}^k\|_{L^2(I_T,H^2)}, \\
&\|( {\ml{A}}^k- {\ml{A}}^{k-1})_t  \|_i \lesssim_0 \|\nabla \bar{u}^k\|_i + T^{1/2}
 \|\nabla(u^k,u^{k-1})\|_1 \| \nabla \bar{u}^k\|_{L^2(I_T,H^2)}\mbox{ for }i=0,\ 1, \\
& \|\nabla \bar{u}^{k}\|_0\leqslant T^{1/2}\| \nabla \bar{u}^k_t\|_{L^2(I_T,L^2)},
\end{aligned}
$$
we make use of \eqref{n0531}, \eqref{n053112} and the above three estimates to deduce from \eqref{difeequion} that
\begin{align}
 \|\nabla \bar{q}^{k+1}\|_{L^\infty(I_T,H^1)} \lesssim_\kappa & T^{1/4}(1+B_2^2) (\|\nabla\bar{u}^k\|_{L^2(I_T,H^2)}\nonumber \\
 &+\| \nabla \bar{u}^k_t\|_{L^2(I_T,L^2)}+ \|\nabla\bar{u}^k\|_{C^0(\overline{I_T},H^1)})\label{202001231907}
\end{align}
and
\begin{align}
    \|\nabla^2  \bar{u}^{k+1}\|_{C^0(\overline{I_T},L^2)}^2+ \|\nabla^3\bar{u}^{k+1} \|_{L^2( {I_T}, L^2)}^2
 \lesssim_\kappa {T} (1+B_2^4)\|\nabla \bar{u}^k \|_{L^2(I_T,H^2)}^2  + \|\nabla  \bar{q}^{k+1}\|_{L^\infty(I_T,H^1)}^2.
 \label{2020012319071}
\end{align}
In addition,
\begin{align}
& \| \bar{u}_t^{k+1}\|_{L^\infty(I_T,L^2)} + \| \bar{u}_t^{k+1}\|_{L^2(I_T,H^1)} \nonumber \\
& \quad \lesssim_\kappa  \|\nabla^2  \bar{u}^{k+1}\|_{C^0(\overline{I_T},L^2)} + (1+B_2^2)
 \| \nabla \bar{u}^k \|_{L^2(I_T,H^2)}  + \|\nabla  \bar{q}^{k+1}\|_{L^\infty(I_T,H^1)} .
\label{202001232104}
\end{align}
%%%
Noting that $( \bar{u}^{k+1})_{\mathbb{T}^3}=0$ and $(\bar{q}^k)_{\mathbb{T}^3}=0$, we put \eqref{202001231907}--\eqref{202001232104}
together to conclude
\begin{equation*}
\begin{aligned}
  \| \bar{u}^{k+1} \|_{\mathfrak{U}_T} + \|\bar{q}^{k+1}\|_{L^\infty(I_T,H^2)}
  \leqslant \tilde{c} T^{1/4}\| \bar{u}^{k+1} \|_{\mathfrak{U}_T}
\end{aligned}
\end{equation*}
for some constant $\tilde{c}$ depending on $B_2$, $\rho$, $\mu$ and $\kappa$.

Finally, we further require $T\leqslant (2\tilde{c})^{-4}$ to find that
$$ \|\bar{u}^{k+1}\|_{\mathfrak{U}_T} +\|\bar{q}^{k+1}\|_{L^\infty(I_T,H^2)}\leqslant\|\bar{u}^{k}\|_{\mathfrak{U}_T}/{2}\quad
\mbox{ for any }\; k\geqslant 1,$$
which implies
\begin{align}
\sum_{k=2}^\infty(\|
\bar{u}^k \|_{\mathfrak{U}_T} +  \| \bar{q}^{k}\|_{L^\infty(I_T,H^2)}  )  <\infty. \nonumber
\end{align}
Hence, $\{ {u}^k,q^k\}_{k=1}^\infty$ is a Cauchy sequence in ${\mathfrak{U}_T}\times L^\infty(\mathbb{R}^+,\underline{H}^2)$ and
\begin{equation}\label{strongconvegneuN}
(\eta^k,u^k, q^k)\rightarrow (\eta ,u, q)  \end{equation}
strongly in $C^0(\overline{I_T},H^3)\times {\mathfrak{U}_T}\times L^\infty(\mathbb{R}^+,\underline{H}^2)$,
where $\eta :=\int_0^tu\mm{d}\tau+\eta^0$.

Consequently, we easily see from \eqref{iteratiequat} and \eqref{strongconvegneuN} that $(\eta,{u},q)$ constructed above
is a solution to the initial value problem \eqref{01dsaf16asdfasf}--\eqref{01dsaf16asdfasfsaf}. The uniqueness of solutions
to \eqref{01dsaf16asdfasf}--\eqref{01dsaf16asdfasfsaf} in the function class
$C^0(\overline{I_T},H^3)\times {\mathfrak{U}_T}\times L^\infty(\mathbb{R}^+, \underline{H}^2)$
can be easily verified by a standard energy method, and its proof will be omitted here.

\vspace{4mm} \noindent\textbf{Acknowledgements.}
The research of Fei Jiang was supported by NSFC (Grant No. 11671086)  and the research of
Song Jiang by NSFC (Grant Nos. 11631008, GZ1465, 11571046)

\renewcommand\refname{References}
\renewenvironment{thebibliography}[1]{ %
\section*{\refname}
\list{{\arabic{enumi}}}{\def\makelabel##1{\hss{##1}}\topsep=0mm
\parsep=0mm
\partopsep=0mm\itemsep=0mm
\labelsep=1ex\itemindent=0mm
\settowidth\labelwidth{\small[#1]}%
\leftmargin\labelwidth \advance\leftmargin\labelsep
\advance\leftmargin -\itemindent
\usecounter{enumi}}\small
\def\newblock{\ }
\sloppy\clubpenalty4000\widowpenalty4000
\sfcode`\.=1000\relax}{\endlist}
\bibliographystyle{model1b-num-names}
 
\end{document}